\documentclass[11pt]{article}

\usepackage[a4paper, total={6in, 8in}]{geometry}

\newcommand{\mc}[1]{\mathcal{#1}}
\newcommand{\bs}[1]{\boldsymbol{#1}}
\newcommand{\cbrac}[1]{\left\{#1\right\}}
\newcommand{\sbrac}[1]{\left[#1\right]}

\usepackage{dsfont}
\newcommand{\indic}[1]{\mathds{1}\brac{#1}}
\usepackage{subfigure}
\usepackage{xcolor}
\usepackage{verbatim}
\usepackage{tikz}
\usepackage{hyperref}
\usepackage{bm}
\usepackage{amsmath}
\usepackage{amsthm}	
\usepackage{amsfonts}	
\usepackage{amssymb, bbm, units}
\usepackage{enumerate}
\usepackage[nobysame, alphabetic, initials]{amsrefs}
\usepackage{ulem} 

\newcommand{\norm}[1]{\left\lVert#1\right\rVert}
\newcommand{\abs}[1]{\left\lvert #1 \right\rvert}
\numberwithin{equation}{section}
\numberwithin{figure}{section}
 \usepackage[nodayofweek]{datetime}

\theoremstyle{definition}
\theoremstyle{plain}
\newtheorem{theorem}{Theorem}[section]
\newtheorem{lemma}[theorem]{Lemma}

\newtheorem{corollary}[theorem]{Corollary}

\newcommand{\EE}{\mathbb{E}}

\newcommand{\PP}{\mathbb{P}}
\newcommand{\RR}{\mathbb{R}}

\newcommand{\brac}[1]{\left(#1\right)}
\newcommand{\mf}[1]{\mathbf{#1}}

\newcommand{\floor}[1]{\left\lfloor{#1}\right\rfloor}
\newcommand{\ceil}[1]{\left\lceil{#1}\right\rceil}
\newcommand\numeq[1]%
  {\stackrel{\scriptscriptstyle(\mkern-1.5mu#1\mkern-1.5mu)}{=}}

\makeatletter
\newsavebox{\@brx}
\newcommand{\llangle}[1][]{\savebox{\@brx}{\(\m@th{#1\langle}\)}%
  \mathopen{\copy\@brx\mkern2mu\kern-0.9\wd\@brx\usebox{\@brx}}}
\newcommand{\rrangle}[1][]{\savebox{\@brx}{\(\m@th{#1\rangle}\)}%
  \mathclose{\copy\@brx\mkern2mu\kern-0.9\wd\@brx\usebox{\@brx}}}
\makeatother

\definecolor{darkgreen}{rgb}{0,0.35,0}

\DeclareMathOperator{\poly}{poly}
\usepackage{hyperref}
\newcommand{\ErdosRenyi}{Erd\H{o}s--R\'enyi }

\DeclareMathOperator{\vol}{Vol}

\newcommand{\ubar}[1]{\underline{#1}}

\title{Robustness of the 2-Choices Dynamics to Node Failures}
\author{Luke Meredith and Arpan Mukhopadhyay\\
University of Warwick}

\date{}

\begin{document}
\maketitle

\begin{abstract}
In many applications, it becomes necessary for a set of distributed network nodes to agree on a common value or {\em opinion} as quickly as possible and with minimal communication overhead.
The classical {\em 2-choices rule} is a well-known distributed algorithm designed to achieve this goal. Under this rule, each node in a network updates its opinion at random instants by sampling two neighbours uniformly at random and then adopting the common opinion held by these neighbours if they agree.
For a {\em sufficiently well-connected network} of $n$ nodes and {\em two} initial opinions, this simple rule results in the network being absorbed in a consensus state in $O(\log n)$ time (with high probability) and the consensus is obtained on the opinion held by the majority of nodes initially.

In this paper, we study the robustness of this  algorithm to node failures. In particular, we assume that with a constant probability $\alpha$, a node may fail to update according to the 2-choices rule and erroneously adopt any one of the two opinions uniformly at random. This is a strong form of failure under which the network can no longer be absorbed in a consensus state. However, we show that as long as the error probability $\alpha$ is less than a threshold value, the network is able to retain the majority support of the initially prevailing opinion for an exponentially long time ($\Omega(\poly(\exp(n)))$). In contrast, when the error probability is above a threshold value, we show that any  opinion quickly ($O(\log n)$ time) loses its majority support and the network reaches a state where  (nearly) an equal proportion of nodes support each opinion. We establish the above phase transition in the dynamics for both complete graphs and expander graphs with sufficiently large spectral gaps and sufficiently homogeneous degrees.  
Our analysis combines spectral graph theory with Markov chain mixing and hitting time analyses.  

\end{abstract}

\section{Introduction}

Collective decision-making models are used throughout statistical physics~\cite{redner_two_choices,chen2005majority} and computer science~\cite{doerr2011stabilizing,cruciani2021phase} to understand how local interaction rules can affect consensus formation. One main focus in this area is the robustness of these collective decision making processes to random node failures. In wireless sensor networks, sensor nodes are prone to anomalous behaviour~\cite{tay2008impact}. In biological networks, agents can deviate randomly from the norm~\cite{briat2023noise}. The key question here is: {\em can we still achieve consensus or consensus-like behaviour in a network when nodes are prone to random failures?} Understanding this fault tolerance allows us to design robust distributed protocols that can achieve desired global objectives. 

In this paper, we study the behaviour of a well-known distributed protocol, namely the \textit{2-choices rule}, for achieving consensus under node failures. Under the classical 2-choices rule~\cite{redner_two_choices,cooper_two_choices}, a node samples two nodes from its neighbourhood uniformly at random and updates to the majority opinion among the sampled nodes and itself. This simple algorithm has many desirable properties. First of all,
it ensures that the system is absorbed in a consensus state very quickly (in $O(\log n )$ time for a network of $n$ nodes) if the nodes are sufficiently well connected. Secondly, consensus is achieved with high probability on the opinion having the initial majority support. This is particularly desirable when the nodes are trying to compute the majority opinion in a distributed way~\cite{vojnovic}.
Furthermore, for graphs with large neighbourhood sizes, this algorithm vastly reduces communication overhead associated with computing the majority opinion, since a node does not need to store the opinions of all nodes in its neighbourhood~\cite{mukhopadhyay2024phase}.  However, many of these desired properties of the 2-choices rule depend on the assumption that nodes behave ideally without failure. In the past, there have been attempts to characterise the effects of failures by introducing a weak adversary capable of changing the opinions of at most $O(\sqrt n)$ number of nodes arbitrarily per unit time~\cite{doerr2011stabilizing,ghaffari2018nearly}. Under such weak adversarial models, no significant change in the behaviour of the 2-choices dynamics have been observed on complete graphs~\cite{doerr2011stabilizing}.

In this paper, we introduce a stochastic model of failure where each node can {\em fail} to update according to the 2-choices rule with a constant probability $\alpha \in (0,1)$ and, in the case of failure, it simply adopts one of the opinions with equal probability. Hence, in our model, $\Theta(n)$ nodes fail per unit of time as opposed to only $O(\sqrt n)$ node failures assumed in earlier works. Our objective is to investigate if the 2-choices model is still robust to such failures. Note that under this model of failure, the network can no longer be absorbed into a consensus state.
Thus, as a measure of robustness, we consider how long the network is able to retain a state close to its initial state, i.e., the time for which the opinion holding the initial majority support continues to do so. We choose this as our measure of robustness since in many real-life situations perfect consensus is not achievable. Instead, achieving majority support for a single opinion and maintaining that majority for a long duration is considered sufficient~\cite{cruciani2021phase}.

Thus, the metrics for robustness we use in this paper are (i) the first time the network hits the state where equal proportions of nodes support each opinion and (ii) the mixing time of the dynamics, which measures how quickly the dynamics reaches its equilibrium distribution. 
The larger these metrics are, the longer the network is able to retain its closeness to the initial distribution where one opinion has a clear majority. 
Our key finding is that the 2-choices dynamics, under the  model of node failure described above, is able to maintain its closeness to the initial state for a time that scales exponentially with the size of the network $n$, as long as the failure probability is less than a constant threshold value. However, above this threshold value, we show that the initial majority information is lost within a time that scales logarithmically with the size of the network $n$. Importantly, we establish this sharp {\em phase transition} in the dynamics not only for {\em dense graphs}, where the average degree of a node scales with the network size, but also for {\em sparse graphs} where the average degree remains bounded as the network scales. Thus, our main result is that the 2-choices model is robust to node failures for a large class of sufficiently well connected networks as long as the failure probability is smaller than a threshold value. To the best of our knowledge, our work is the first to show such strong robustness guarantees for the 2-choices dynamics for both dense and sparse graphs.  Here are our key contributions in more detail.
\begin{itemize}
\item We start by analysing the  2-choices dynamics with node failures on complete graphs. Here, we show that as long as the failure probability $\alpha$ is less than $1/3$, it takes an exponentially long time to disrupt the initial majority and the mixing time of the dynamics also scales exponentially with $n$. However, for $\alpha > 1/3$ we show that both of these times scale logarithmically with the size of the network. 

\item Next, we generalise the above phase transition result to a large class of graphs with sufficiently large spectral radius and sufficiently homogeneous degree distributions, which includes both dense and sparse  graphs. In particular, for dense graphs with degrees growing with $n$, our results show that the threshold for the phase transition remains the same as complete graphs (1/3). Furthermore, for sparse graphs with bounded degrees, we find a non-zero threshold value for the error probability $\alpha$, below which, the network is guaranteed to be robust to failures. 

\item To establish our results, we prove several general lemmas that are helpful in analysing the mixing and hitting time of continuous-time birth-death processes. Furthermore, using spectral methods, we prove tight bounds on transition rates of the dynamics on general graphs. Our approach is sufficiently general to be useful in analysing similar non-linear dynamical systems on complex networks.
\end{itemize}

\subsection{Related Literature}

The earliest model in distributed systems designed to achieve consensus is the voter model \cite{clifford1973model, holley1975ergodic}, where a node copies the opinion of a sampled neighbour. Therefore, the probability a node switches opinion is equal to the proportion of nodes in its neighbourhood holding the opposing opinion. Due to the simplicity of these dynamics, the Voter model has been extensively studied and applied in a variety of areas such as statistical physics and social sciences \cite{redner2019reality}. The model has been studied on different classes of graphs including lattices \cite{cox1989coalescing}, random $d$-regular graphs \cite{cooper2013coalescing}, and \ErdosRenyi graphs \cite{nakata1999probabilistic}, due to the duality between the linearity of the model and coalescing random walks \cite{holley1975ergodic, clifford1973model}. It is a known result that the Voter model reaches consensus in the initial majority in $O(n)$ time on expander graphs, and the probability of achieving this consensus is small \cite{cooper2015fast}. Noise was first introduced into the Voter model in \cite{noisy_voter_model}, where nodes spontaneously change state, independent of the regular update rule, and critical exponents associated with a phase transition were derived. 

Another cornerstone in distributed systems is the majority-rule model, where the updating node takes the majority opinion within its entire neighbourhood. Hence, this model is designed to drive the system towards its initial majority. A generalization of this model involves a node sampling $m$ neighbours and switching opinion if $d$ or more of them differ in opinion from its own opinion. This was first studied in \cite{cruise2014probabilistic}. Thus, the 2-choices rule is a simplification of this model. The particular model of the classic 2-choices rule used in this paper was first analyzed on random regular graphs and expanders in \cite{cooper_two_choices,cooper2015fast}. Under this model, it has been shown that the network reaches consensus on the opinion having the initial majority with high probability, in $O(\log n)$ time, as long as the initial proportions of nodes supporting the two opinions are sufficiently different. The 2-choices model under an $F$-bounded adversary was studied in~\cite{doerr2011stabilizing}. The $F$-bounded adversary is assumed to be capable of arbitrarily changing at most $F$ opinions at each discrete time step. It was shown that for $F=O(\sqrt{n})$, the two choices dynamics can still achieve consensus on the initial majority opinion in $O(\log n)$ time for complete graphs. 

Many recent papers~\cite{vieira2016phase, cruciani2021phase, mukhopadhyay2024phase, d2022phase} study similar dynamical models for biased models of noise, where one opinion is favoured over another. In \cite{vieira2016phase}, a model with a regular update mechanism and an independent noise is studied. Our model uses this fundamental framework, with the regular update mechanism being the classical 2-choices rule. Our key focus though is to identify regimes where the system is robust to this independent noise on a class of complex networks. Our use of both the hitting time to break the initial majority and the mixing time as measures of robustness against noise is one aspect of where these two approaches differ. In \cite{cruciani2021phase}, a `Biased-2-Choices($p, \sigma$) dynamics' is introduced whereby a node's sampled neighbour supports the colour/opinion $\sigma$ of the node with probability $p$. This type of bias can be viewed as a communication noise between two nodes, and is one directional, i.e. it favours one opinion over another. In \cite{mukhopadhyay2024phase}, we see a similar noise model, where one opinion is viewed superior to another using a bias. In this model, a node will spontaneously change state without comparing opinions with other players. This mechanism is closer to the one studied in this paper, but the addition of bi-directional noise in our work renders a departure from this analysis. In \cite{d2022phase}, the authors study a closely related model, where the noise is embedded in transmission. This model is a 3-choice majority model, whereby a node samples three neighbours who support an opinion that gets corrupted with noise. The initially sampled node then takes the majority in the group. Interestingly, the authors obtain the same critical value where the phase transition occurs on complete graphs as we do for the model in this paper. However, in their model, one opinion is preferred over another. Furthermore, their analysis does not extend to sparse graphs where the degrees of nodes can be bounded.

Despite extensive work on distributed consensus protocols, our paper is the first to study the 2-choices dynamics under a strong model of failure where $\Theta(n)$ nodes can fail in each unit time. Also, unlike most previous papers which study noisy opinion dynamical systems, we are the first to prove results that hold for both sparse and dense graphs.


\subsection{Organisation}

The remaining sections of the paper are organised as follows. In Section \ref{Section: Model}, we define the 2-choices model with node failures studied in the paper. In Section \ref{Section: Measures of Robustness}, we characterise the stationary distribution of the model to argue why hitting times and mixing times are good measures of robustness. Next, in Section \ref{Section: Complete graphs}, we study the dynamics on complete graphs, and in Section \ref{Section: General graphs}, we extend this to general graphs. In Section \ref{sec:gen_proofs}, we provide proofs of the general lemmas  used in the proofs.  Section \ref{Section: Conclusion} provides a summary of the results and concludes the paper.

\section{Model}
\label{Section: Model}

We consider a network of $n$ nodes situated on the vertices of a simple,  undirected, and connected graph $G=(V,E)$ where $V$, with $\abs{V}=n$, denotes the set of vertices and $E$ denotes the set of edges. For two connected nodes $u,v\in V$ the edge between them is denoted by the unordered pair $\{u,v\}$. The set of neighbours of a node $v$ is denoted by $N_v=\{u\in V: \{u,v\}\in E\}$. 

At any time $t\geq 0$, each node $v\in V$ is assumed to have an opinion $X_v(t)$ in the set $\{0,1\}$. The opinions of nodes evolve as follows.
Each node updates its opinion at points of an independent unit rate Poisson process associated with itself. At an instant of update, one of the two possible events occurs:  (i) with probability $1-\alpha$ the updating node samples two of its neighbours uniformly at random (with replacement), and only if these sampled neighbours are found to have the same opinion, the updating node adopts the common opinion of the sampled neighbours (otherwise, the updating node retains its previous opinion), (ii) with probability $\alpha$ independently of any nodes' opinion (including itself) the updating node adopts one of the two opinions uniformly at random. In the former case, we refer to the update event as a {\em 2-choices update} and in the later case we refer to the update event as a {\em node failure}. Thus, the parameter $\alpha$ denotes how often a node {\em fails} to update according to the 2-choices rule. 

We denote the random state of the network at time $t \geq 0$ by $\mf X(t)=(X_v(t), v \in V)\in \mc{X}:=\{0,1\}^n$. Clearly, the process $\mf X=(\mf X(t) , t \geq 0)$ is Markov and takes values in $\mc{X}$. We denote specific realisations of the random state vector by smaller case letters, e.g., $\mf x, \mf y$ etc. We also denote the set of all possible subsets of the state space $\mc X$ by $\mc{P}(\mc X)$.
In a given state $\mf x \in \mc{X}$, the set of nodes with opinion $1$ will be denoted by $A(\mf x)=\{v\in V: x_v=1\}$. Similarly, $B(\mf x)=V-A(\mf x)$ will denote  the set of nodes with opinion $0$ in state $\mf x$. 
With a slight abuse of notation, we will also use $A(\mf x)$ and $B(\mf x)$ to denote the sizes of the sets $A(\mf x)$ and $B(\mf x)$, respectively. 
Clearly, the process $\mf X$ transitions from state $\mf x$ to state $\mf y\neq \mf x$ with a rate given by

\begin{equation}
    q(\mf x,\mf y)=\begin{cases}
        \brac{\frac{\alpha}{2}+(1-\alpha)\brac{\frac{d_v^{A}(\mf x)}{d_v}}^2}\indic{x_v=0}, &\text{if } \mf y =\mf x +\mf e_v \text{ for some } v \in V,\\
        \brac{\frac{\alpha}{2}+(1-\alpha)\brac{\frac{d_v^{B}(\mf x)}{d_v}}^2}\indic{x_v=1}, &\text{if } \mf y =\mf x -\mf e_v \text{ for some } v \in V,\\
        0, &\text{otherwise,}
    \end{cases}
    \label{eq:rates_graph}
\end{equation}
where $\mf e_v$ denotes the $n$ dimensional unit vector with $1$ in the $v^{\textrm{th}}$ position, $d_v=\abs{N_v}$ denotes the degree of node $v$, $d_v^S=\abs{N_v \cap S}$ denotes the set of neighbours of $v$ in $S$ for any $S\subseteq V$, and $d_v^A(\mf x)$ and $d_v^B(\mf x)$ denote $d_v^S$ for $S=A(\mf x)$ and $S=B(\mf x)$, respectively. 

The case $\alpha=0$ corresponds to the classical 2-choices model where there are two absorbing states - the all-$1$ state, denoted by $\bm{1}$, and the all-$0$ state, denoted by $\bm{0}$. It is well-known~\cite{cooper_two_choices,cooper2015fast} that the classical 2-choices dynamics converge to the absorbing state closest (in $\ell^1$ norm) to the initial state $\mf X(0)$ in $O(\log n)$ time (with high probability) if the underlying graph $G$ is {\em sufficiently well-connected} (more precisely, if $G$ is an expander graph) and the proportions of nodes supporting each opinion are sufficiently well-separated at $t=0$. Hence, in the classical 2-choices dynamics, the initial majority opinion continues to hold its majority support at all times, eventually becoming the opinion at {\em consensus}.
However, when $\alpha >0$, it is easy to see that there are no absorbing states in the dynamics.
Indeed, the process $\mf X$ is ergodic  with a unique equilibrium distribution which we denote by $\bs \pi$.
Although the process $\mf X$ seems to  behave very differently based on whether $\alpha > 0$ or $\alpha=0$, we shall try to characterise how similar (resp. different) the dynamics is to (resp. from) the classical 2-choices dynamics for smaller (resp. larger) values of $\alpha$.
 
{\bf General Notations}: We use the following notations throughout the rest of the paper. We use $\PP_{\mf x}(\cdot)$ and $\EE_{\mf x}[\cdot]$ to respectively denote the probability and expectation conditioned on $\mf X(0)=\mf x$. The total variation distance between any two distributions $\mu$ and $\nu$ on the measurable space $(\mc{X}, \mc{P}(\mc{X}))$ will be denoted by $\norm{\mu-\nu}_{\text{TV}}:=\sup_{\mc H\subseteq \mc{P}(\mc{X})}\abs{\mu(\mc H)-\nu(\mc H)}$.  For $E_1,E_2 \in \mc{P}(\mc{X})$, we say that {\em $E_1$ occurs with high probability given $E_2$} when $\PP(E_1|E_2)\to 1$ as $n\to \infty$.  
For any $x,y, z\in \RR$ we let $x\vee y=\max(x,y)$, $x\wedge y=\min(x,y)$, $[z]_+=z\vee 0$ and $[z]_-=-z\vee0$. We write $f(n)=O(g(n))$ if $\lim\sup_{n\to\infty}(f(n)/g(n)) < \infty$; $f(n)=o(g(n))$ if $\lim_{n \to \infty}f(n)/g(n) = 0$; $f(n)=\Omega(g(n))$ if $g(n)=O(f(n))$; and $f(n)=\Theta(g(n))$ if $f(n)=O(g(n))$ and $f(n)=\Omega(g(n))$. We use $\indic{\cdot}$ to denote the indicator function. Furthermore, for two random variables $X$ and $Y$ we write $X\geq_{{st}} Y$ if $\PP(X\geq x)\geq \PP(Y\geq x)$ for all $x\in \RR$.

\section{Measures of Robustness}
\label{Section: Measures of Robustness}

In this section, we discuss how we measure the robustness of the 2-choices dynamics to node failures. To define and justify the measures of robustness, we first need to characterise the stationary distribution $\bs{\pi}$ of the Markov process $\mf X$ defined in the previous section. An important observation that follows from the symmetry of the dynamics is given in the following lemma.
\begin{lemma}
\label{lem:stationary}
For any $\mf x\in \{0,1\}^n$, we have $\bs \pi(\{\mf x\})=\bs \pi(\{\bm{1}  -\mf x\})$. Hence, $\EE_{\mf X \sim \bs \pi}[A(\mf X)]=\frac{n}{2}$.
\end{lemma}

Thus, for all positive values of the failure probability $\alpha$, the stationary distribution of $A(\mf{X})$ is symmetric around $n/2$.  This implies that, as long as there are node failures, the states where $1$ is the majority opinion are visited equally often as the states where $0$ is the majority opinion in the long-run.
This is in sharp contrast to the classical 2-choices dynamics without node failures (i.e., $\alpha=0$) where the majority opinion at $t=0$ continues to be the majority opinion at all times, eventually becoming the opinion at consensus, i.e., the dynamics is always contained within one side of $n/2$. 
Thus, a natural way to measure the robustness of the 2-choices dynamics to node failures is through  the time the process spends on one side of $n/2$ before switching to the other side for the first time. More precisely, if $T_p:=\inf\cbrac{t\geq 0: A(\mf X(t))=\ceil{pn} \text{ or } A(\mf X(t))=\floor{pn}}$ for $p\in(0,1)$ denotes the first {\em hitting time} to a state where (approximately) $p$ fraction of nodes support opinion $1$, then we can consider $T_{1/2}$ as a measure of robustness to failures. The higher the value of $T_{1/2}$, the longer the network is able to retain the majority of the initially prevailing opinion and more similar the dynamics is to the classical 2-choices dynamics (making it more robust to node failures).

Another way to measure closeness to the initial distribution is through the {\em mixing time} which gives the rate at which the dynamics converges to its equilibrium distribution. 
Recall that, for a fixed $\epsilon >0$, the $\epsilon$-mixing time $t_{\text{mix}}(\epsilon)$ is defined as
$t_{\text{mix}}(\epsilon)=\inf\cbrac{t \geq 0:d(t) \leq \epsilon}$,
where $d(t)=\max_{\mf x \in \mc{X}}\max_{\Gamma \in \mc{P}(\mc X)}\abs{\PP_{\mf x}(\mf X(t)\in \Gamma)-\bs{\pi}(\Gamma)}$ is the maximum possible total variation distance between the distribution of $\mf X(t)$ and the stationary distribution $\bs \pi$ where the maximum is taken over all possible choices of the initial state~\cite{peres_book}.
%
Thus, $t_{\text{mix}}(\epsilon)$ denotes the first time the total variation distance between the distribution of $X(t)$ and $\bs \pi$ falls below $\epsilon$. As is standard in the literature~\cite{peres_book}, we define the mixing time by $t_{\text{mix}}:=t_{\text{mix}}\brac{\frac{1}{4}}$.  The higher the mixing time, the longer the network remains close to its initial distribution resulting in a network more robust to node failures.

The following corollary, which is a consequence of Lemma~\ref{lem:stationary}, will be useful in connecting a lower bound on $T_{1/2}$ to a lower bound on $t_{\text{mix}}$.

\begin{corollary}
\label{cor: TV distance bound}
For any $\mf x\in \mc{X}$, we have
\begin{equation}
d(t)\geq \frac{1}{2}-\PP_\mf{x}(T_{1/2}\leq t).
\end{equation}
\end{corollary}

\begin{proof}
By definition of total variation distance and $d(t)$, for any $\mf x \in \mc{X}$ and any $\mc H\subseteq \mc X$ we have
\begin{align*}
d(t) 
&\geq \abs{\PP_{\mf x}(\mf X(t)\in \mc{H})-\bs{\pi}(\mc{H})}.
\end{align*}
%
Now, we choose the set $\mc H\equiv\mc H(\mf x)$ depending on $\mf x$ as follows: if $A(\mf x)\geq \ceil{n/2}$ we choose $\mc H(\mf x)=\{\mf y \in \mc X: A(\mf y)\leq \ceil{n/2}\}$ and if $A(\mf x)\leq \floor{n/2}$, then we choose $\mc H(\mf x)=\{\mf y\in \mc X: A(\mf y)\geq \floor{n/2}\}$. Hence, we always have $\bs \pi(\mc H(\mf x))\geq 1/2$ according to Lemma~\ref{lem:stationary}. This implies that for any $\mf x\in \mc X$ we have
\begin{align}
d(t)&\geq \abs{\PP_{\mf x}(\mf X(t)\in \mc{H}(\mf x))-\bs{\pi}(\mc{H}(\mf x))}\\
&\geq (\bs{\pi}(\mc{H}(\mf x))-\PP_{\mf x}(\mf X(t)\in \mc{H}(\mf x)))\\
&\geq \frac{1}{2}-\PP_{\mf x}(T_{1/2}\leq t),
\end{align}
where the third line follows since $\abs{u}\geq  -u\vee u$ for any $u\in \RR$ and the last line follows from the definition of $T_{1/2}$.
\end{proof}

\section{Complete Graphs}
\label{Section: Complete graphs}

In this section, we analyse the dynamics on complete graphs where each node is connected to every other node in the network, i.e.,  we have $N_v=V-\{v\}$ for all $v\in V$. In this case, the dynamics can be fully described in terms of only the total number of nodes with opinion $1$ at each instant, i.e., the process $A=(A(t), t\geq 0)$ where $A(t):= A(\mf X(t))$ is itself Markov, taking values in $\{0,1,\ldots,n\}$. It follows easily that $A$ is a birth-death process with birth and death rates in state $a\in \{0,1,\ldots,n\}$ given by 
\begin{align}
    q_+(a) &= (n-a) \left[ \frac{\alpha}{2} + (1-\alpha)\left( \frac{a}{n-1}\right)^2 \right],\label{eq:up}\quad \text{and}\\
    q_-(a) &= {a} \left[ \frac{\alpha}{2} + (1-\alpha)\left( \frac{n-a}{n-1}\right)^2\right],\label{eq:down}
\end{align}
respectively. Although this reduction in the state-space dimension of the underlying Markov process considerably simplifies our analysis, as we shall see in the next section, the key insights we obtain from this analysis also hold for general graphs. Thus, the analysis presented in this section for complete graphs forms the basis for the more general analysis presented in later sections. 

Our main results for complete graphs are summarised in the following theorem.

\begin{figure}
    \centering
    \includegraphics[width=\linewidth]{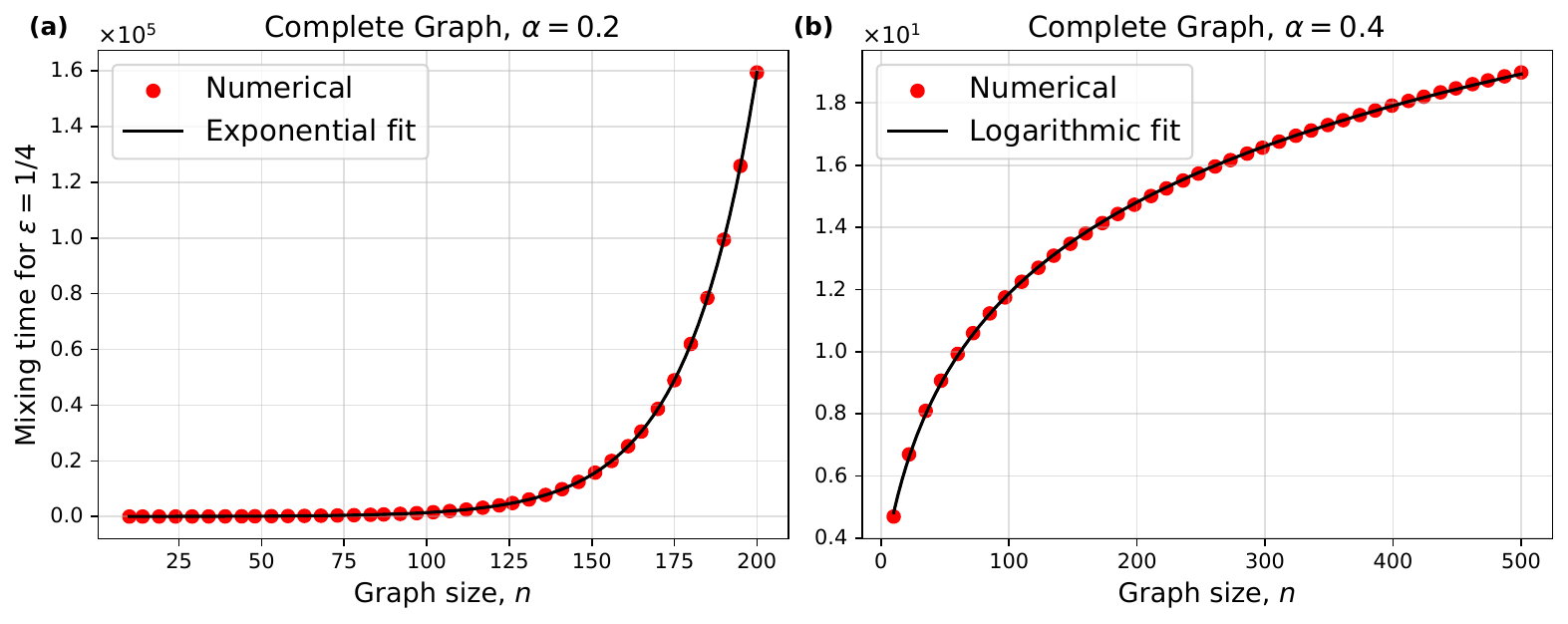}
    \caption{Mixing time ($t_{\text{mix}}$) as a function of the network size ($n$) for different values of the failure probability ($\alpha$) for complete graphs.  \textbf{(a)} $\alpha=0.2 <1/3$. \textbf{(b)} $\alpha=0.4 > 1/3$. In both cases, to compute $t_{\text{mix}}$, the transient distribution of $A(t)$  is found by solving the forward Kolmogorov equation, and the stationary distribution $\bs{\pi}$ is calculated by solving the detailed balance equations $\bs{\pi}(\{a\})q_+(a)=\bs{\pi}(\{a+1\})q_-(a+1)$ for $a\in \{0,1,\ldots,n-1\}$.}
    \label{fig:mixing times complete}
\end{figure}

\begin{theorem}
\label{thm:complete}
For the 2-choices dynamics with node failures on complete graphs, the following properties hold.
\begin{enumerate}
\item For $0 < \alpha < 1/3$, we have $t_{\text{mix}}= \Omega(\exp(\Theta(n)))$ and $T_{1/2}=\Omega(\exp(\Theta(n)))$ with high probability given that $A(0)=\ceil{np}$ for some $p\neq 1/2$.

\item For $\alpha > 1/3$, we have $t_{\text{mix}}=O(\log n)$ and $\EE_{a}[T_{1/2}]=O\mc (\log n)$ for any $a\in \{0,1,\ldots,n\}$ where $\EE_{a}[\cdot]$ denotes expectation conditioned on the event $A(0)=a$. 
\end{enumerate}
\end{theorem}
{\em Discussion on the main result}: According to the theorem above,  the dynamics of the 2-choices model with node failures undergo a {\em phase transition} at $\alpha=1/3$.
When the failure probability $\alpha < 1/3$, both the mixing time $t_{\text{mix}}$ and the hitting time $T_{1/2}$ (starting from any state $\ceil{np}$ with $p\neq 1/2$) scale exponentially with the network size $n$,
implying that the network spends an exponentially long time (in $n$)
on either side of (below or above) the state $n/2$ crossing the $n/2$ boundary only very rarely. This behaviour for smaller values of the failure probability $\alpha$, is  similar to the classical 2-choices model (with $\alpha=0$) where the trajectory of the process never crosses the state $n/2$. Thus, although in the steady-state, the process $A$ spends an equal proportion of time on each side of the state $n/2$ (which follows from the symmetry of $\bs \pi$ established in Lemma~\ref{lem:stationary}), 
when observed on relatively shorter (but still exponentially long) time scales, the process seems to be contained only in one side of $n/2$.
This dichotomy in the dynamics across different time scales is often referred to as {\em metastability} in the literature~\cite{cruciani2021phase,arpan_JSP_20}. 
In sharp contrast, when $\alpha > 1/3$, the process $A$ crosses the $n/2$ state for the first time in only $O(\log n)$ time. This implies that the initial majority is disrupted rapidly. Furthermore, the crossing of the $n/2$ state occurs once every $O(\log n)$ time, giving a very fast convergence (and a small mixing time) to the stationary distribution $\bs{\pi}$. In Figure~\ref{fig:mixing times complete}, we plot the mixing time $t_{\text{mix}}$ of the process $A$ as a function of $n$ for different values of $\alpha$. The results confirm the theory and show that the mixing time is indeed exponential for $\alpha < 1/3$ and logarithmic for $\alpha > 1/3$.

\begin{figure}
    \centering
    \includegraphics[width=\linewidth]{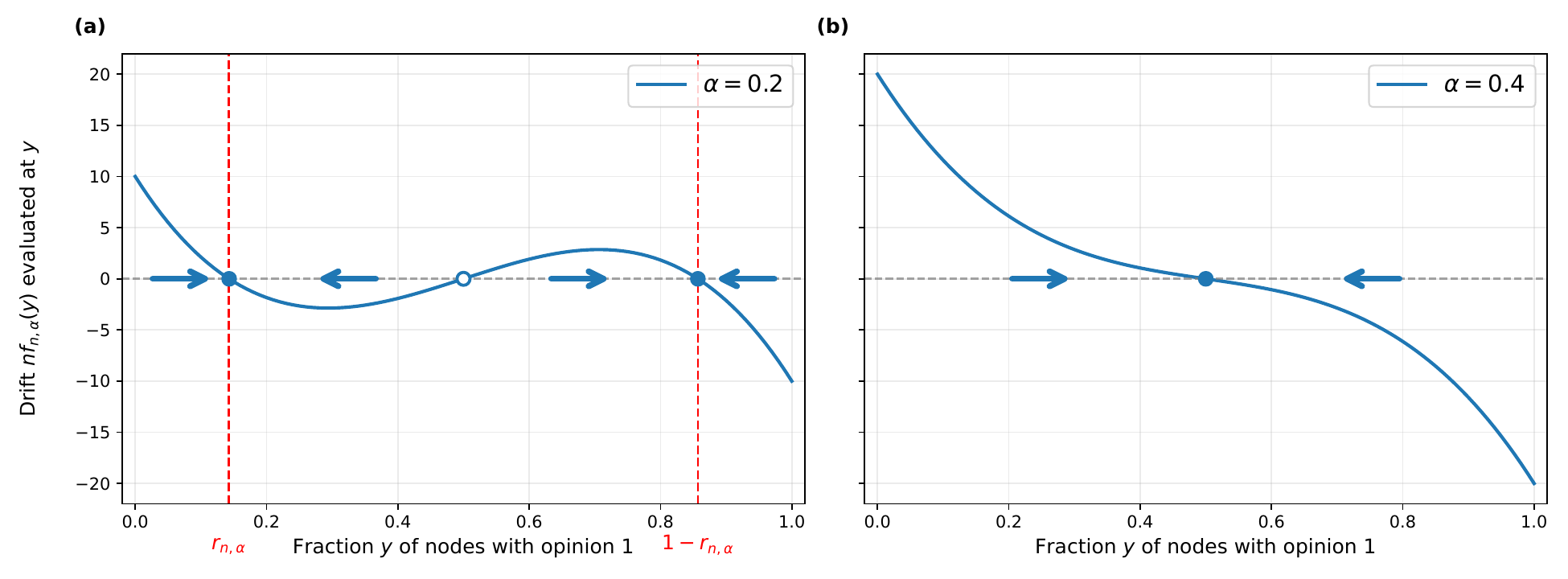}
    \caption{Plot of $f_{n,\alpha}$ as a function of the fraction of nodes with opinion $1$ for different values of $\alpha$.  The arrows in both plots indicate the direction of the drift (right for positive drift and left for negative drift). In both plots, we have chosen $n=100$.  \textbf{(a)} For $\alpha=0.2 <1/3$ there are three distinct real roots of $f_{n,\alpha}$, which, in the increasing order of their values, are at $r_{n,\alpha}, 1/2,$ and $1-r_{n,\alpha}$. All trajectories are attracted either toward $r_{n,\alpha}$ or toward $1-r_{n,\alpha}$.\textbf{(b)} For $\alpha = 0.4>1/3$ the unique root of $f_{n,\alpha}$ is at $1/2$ and all trajectories are attracted toward this root.}
    \label{fig:drift plots}
\end{figure}

{\em Intuitive explanation of the main result using drift}: To explain the above phase transition in the dynamics of $A$, we need to consider {\em the drift} of any function $h:\{0,1,\ldots,n\}\to \RR$ along the trajectory of $A$, which characterises the rate at which $h(A(t))$ changes with $t$. This is defined by the 
generator operator $\mc{G}_A$ of the Markov chain $A$  acting on the function $h$. More precisely, we have
$$\mc{G}_A h(a):=\lim_{h \downarrow 0} \frac{\EE[A(t+h)-A(t)|A(t)=a]}{h}=q_+(a)(h(a+1)-h(a))+q_-(a)(h(a-1)-h(a)),$$ where $q_+$ and $q_-$ are as defined in~\eqref{eq:up}-\eqref{eq:down}. In particular, we are interested in the drift of $h$ when $h=I$, where $I:\RR\to \RR$ is the identity map $a\mapsto I(a)=a$. The drift of $I$ captures the rate at which number of nodes with opinion 1 changes in the network. If this rate is positive (resp. negative) at some state, it implies that the number of nodes with opinion $1$ tends to increase (resp. decrease) on average in that state. Inserting $h=I$ into the definition of the drift we obtain
\begin{align}
    \mc{G}_A I(a)&=q_+(a)-q_-(a)\\
    &= nf_{n,\alpha}\brac{\frac{a}{n}},
\end{align}
where $f_{n,\alpha}:(0,1)\to (0,\infty)$ is defined as
\begin{equation}
    f_{n,\alpha}(y)=(1-2y)\frac{\alpha}{2}-(1-\alpha)\brac{\frac{n}{n-1}}^2y(1-y)(1-2y).
\end{equation}
Now, using the properties of the function $f_{n,\alpha}$ established later in Lemma~\ref{lem:f_complete}, we see that for each $0 < \alpha < 1/3$ there exists $r_{\alpha} \in (0,1/2)$ such that $\mc G_A I(a) < 0$ when $a/n \in (r_\alpha, 1/2)$, and $\mc G_A I(a) > 0$ when $a/n \in (1/2, 1-r_\alpha)$. Hence, whenever the proportion of nodes with opinion $1$ is in the range $(r_\alpha, 1/2)$ (resp. $(1/2, 1-r_\alpha)$), this proportion tends to decrease (resp. increase) on average, making it more difficult to reach the state $n/2$. This drift away from the state $n/2$ causes the hitting time to the state $n/2$ to be very large for $\alpha < 1/3$. In contrast, for $\alpha > 1/3$, Lemma~\ref{lem:f_complete} shows that $\mc G_A I(a) > 0$ for $a/n \in (0,1/2)$ and $\mc G_A I(a) < 0$ for $a/n \in (1/2,1)$. Hence, in this case, the state $n/2$ acts as a global attractor where all trajectories tend to converge quickly. The behaviour of the function $f_{n,\alpha}$ for different values of $\alpha$ is shown in Figure~\ref{fig:drift plots}. It clearly shows that the drift points away from $n/2$ for $\alpha < 1/3$  and toward $n/2$ for $\alpha > 1/3$. 

\textit{Sketch Proof of Theorem \ref{thm:complete}:} 
To prove Theorem~\ref{thm:complete}, we need three main lemmas. 
The first lemma stated below is specific to the Markov chain $A$ and characterises the properties of its drift function $f_{n,\alpha}$, as discussed in the preceding paragraph. 
In particular, it shows that, for $\alpha < 1/3$, $f_{n,\alpha}(y) < 0$ (resp. $f_{n,\alpha}(y) > 0$) whenever $y \in (r_{\alpha},1/2)$ (resp. $y \in (1/2, 1-r_{\alpha})$) where $r_{\alpha} \in (0,1/2)$ is a constant depending on $\alpha$.  Similarly, for $\alpha > 1/3$ and $n$ sufficiently large, the lemma shows that $f_{n,\alpha}(y) > 0$ (resp. $f_{n,\alpha}(y) < 0$) whenever $y\in (0,1/2)$ (resp. $y\in (1/2,1)$). The lemma also shows that, for $\alpha >1/3$ and $n$ sufficiently large, the current proportion $y$ of nodes with opinion $1$ in the network approaches $1/2$  at a rate of at least $c(n,\alpha)\abs{y-1/2}$ for some positive constant $c(n,\alpha)$ depending on both $\alpha$ and $n$. This follows from the first and last statements of the lemma below, since by these we have $f_{n,\alpha}(y)=f_{n,\alpha}(y)-f_{n,\alpha}(1/2)\leq -c(n,\alpha)(y-1/2)$. This exponential rate of convergence is responsible for the fast mixing of the chain in the regime $\alpha >1/3$.

\begin{lemma}
\label{lem:f_complete}
The following statements hold.
\begin{enumerate}
    \item For all $\alpha \in (0,1)$ and all $y\in (0,1)$, we have $f_{n,\alpha}(y)=-f_{n,\alpha}(1-y)$. In particular, we have $f_{n,\alpha}(1/2)=0$.
    \item Fix $\alpha < 1/3$ and let $r_{\alpha}=\frac{1}{2}\brac{1-\sqrt{\frac{1-3\alpha}{1-\alpha}}}$. Then, $f_{n,\alpha}(y) < 0$ when $y\in (r_{\alpha},1/2)$. 
    \item For all $\alpha \in (0,1)$, we have $f_{n,\alpha}(y)-f_{n,\alpha}(y')\leq -c(n,\alpha) (y-y')$ for any $0\leq y'\leq y\leq 1$, where $c(n,\alpha) = (3\alpha-1)/2-o(1)$.  In particular, for $\alpha > 1/3$ and $n$ sufficiently large, we have $f_{n,\alpha}(y) > 0$ for all $y \in (0,1/2)$.
\end{enumerate}
\end{lemma}

To use the above properties of the drift in deriving the bounds on the hitting and mixing time stated in Theorem~\ref{thm:complete}, we need 
the next two lemmas - which will also be useful in the next section.

The following general result for continuous-time birth-death chains defined on the state-space $\{0,1,\ldots,n\}$ states that, within the state-space of the chain, if there is a significantly large region ranging from $\ceil{nr_1}$ to $\floor{nr_2}$ (for some $0< r_1<r_2< 1$) where the drift of the chain is negative (given by condition~\eqref{eq:cond_neg_drft} of the lemma), then with high probability, the chain's hitting time to the state $nr_2$ will be exponentially long in $n$ (as stated in~\eqref{eq:exp_long}),
provided that the holding time in each state is uniformly bounded below by $c/n$ where $c>1$ (given by the condition~\eqref{eq:cond_lazyness}).
The intuition for this result is that, due to its negative drift, the chain visits the state $\ceil{nr_1}$ exponentially many times before it reaches the state $\floor{nr_2}$, spending a sufficiently long time in each such visit. To apply this lemma to the birth-death chain $A$, we will verify using Lemma~\ref{lem:f_complete} and the rates in~\eqref{eq:up}-\eqref{eq:down}, that conditions~\eqref{eq:cond_neg_drft} and~\eqref{eq:cond_lazyness} indeed hold for $A$ with $r_1=r_\alpha$ and $r_2=1/2$, as long as $\alpha < 1/3$. This will establish the first statement of Theorem~\ref{thm:complete}.

\begin{lemma}
\label{lem:exp_bound_gen}
    Let $Y=(Y(t) : t\geq0)$ be a continuous-time birth-death chain on the state space $\{0, 1, \ldots, n\}$ with birth and death rates in state $y\in \{0,1,\ldots,n\}$ given by $q_+(y)=nf_n(y/n)$ and $q_-(y)=ng_n(y/n)$, respectively, where $f_n,g_n:[0,1] \rightarrow \RR_+$ are both continuous functions satisfying 
    \begin{equation}
        f_n(z)-g_n(z) < 0, \forall z\in (r_1,r_2), \label{eq:cond_neg_drft}
    \end{equation}
    for some $0 < r_1<r_2< 1$ and $f_n(z),g_n(z)>0$ for all $z\in (0,1)$. Furthermore, assume that there exist some constants $c>1$ and $n_0 \geq1$ such that for all $y\in \{0,1,\ldots,n\}$ and for all $n\geq n_0$, we have
    \begin{equation}
    \label{eq:cond_lazyness}
        0 < q_-(y) + q_+(y) \leq \frac{1}{c}n.
    \end{equation}
    Also, assume that the function $h_n=g_n/f_n$ satisfies $\inf_{n, z\in [r_1+\epsilon,r_2-\epsilon]}h_n(z) >1$ for any $\epsilon >0$ and has a uniformly bounded derivative in the range $[r_1,r_2]$, i.e., $\sup_{n,z\in[r_1,r_2]} \abs{h_n'(z)} < \infty$.
    Define $T_r = \inf\{t\geq0 : Y(t) = \ceil{nr} \text{ or } Y(t)=\floor{nr}\}$. Then, for each $p\in (r_1,r_2)$ and  $n\geq 2\vee n_0$ we have
    \begin{equation}
        \PP_{\ceil{np}}(T_{r_2}\geq \exp(\kappa_1n))\geq 1-\kappa_2n^2\exp(-2n\kappa_1),
        \label{eq:exp_long}
    \end{equation}
    where $\kappa_1 \text{ and } \kappa_2$ are positive constants and $\PP_{y}(\cdot)$ denotes {the probability} conditioned on $Y(0)=y$.
\end{lemma}

The final lemma  required to complete the proof of Theorem~\ref{thm:complete} is stated below.
It characterises the hitting time of a continuous-time Markov chain to states where a non-negative function $\Phi$, defined on the state-space of the chain, becomes zero assuming that $\Phi$ contracts at a rate of at least $c \Phi$ on average (for some positive constant $c>0$) throughout the state-space of the chain. Specifically, the lemma shows that any state where $\Phi=0$ is reached in $O(\log \Phi_{\max})$ time on average, where $\Phi_{\max}$ is the maximum value of $\Phi$ in the state-space. Intuitively, this occurs due to $\Phi$ decaying exponentially fast along the entire trajectory of the chain. 
To apply this lemma to the birth-death chain $A$, we recall that Lemma~\ref{lem:f_complete} already establishes that the distance between any state $a$ and $n/2$ contracts in the same way as stated in the following lemma. Hence, to reach $n/2$, it only takes $O(\log n)$ time on average.

\begin{lemma}
    \label{lem:log_bound_gen}
    Let $(Y(t) : t\geq0)$ be a continuous-time Markov chain with a finite state space $\mc Y$ and generator $\mc G_Y$. Let $\Phi : \mc{\bar Y} \rightarrow\mathbb{R}_+$ be a function defined on $\mc{\bar Y}\supseteq \mc Y$ such that there exists ${\mc Y}_0 \subseteq {\mc  Y}$ where $\Phi(y)=0$ for all $y\in \mc{ Y}_0$ and
    \begin{align}
        \mc G_Y \Phi(y) \leq - c\Phi(y), \quad \forall y \in \mc Y,
    \end{align}
    for some constant $c>0$. Then, the hitting time $\tau=\inf\{t\geq 0: Y(t)\in \mc Y_0\}$ satisfies 
    \begin{align}
        \EE_y[\tau] \leq \frac{1+\log \brac{{\Phi(y)}/{\Phi_{\min}}}}{c},  \quad  \forall y \in \mc Y-\mc{ Y}_0,
    \end{align}
    where $\Phi_{\min}=\min_{y\in \mc Y-\mc { Y}_0}\Phi(y) > 0$.
\end{lemma}

We now provide a detailed proof of Theorem~\ref{thm:complete} proceeding along the lines discussed above. 

{\em Proof of Theorem~\ref{thm:complete}}: To prove the first statement of Theorem \ref{thm:complete}, 
we apply Lemma~\ref{lem:exp_bound_gen} to the birth-death chain $A$. Note from~\eqref{eq:up}-\eqref{eq:down} that $q_+(a)=nf_n(a/n)$ and $q_-(a)=ng_n(a/n)$ with
\begin{align}
    f_n(y)&=(1-y)\brac{\frac{\alpha}{2}+(1-\alpha)\brac{\frac{n}{n-1}}^2 y^2}, \text{ and}\\
    g_n(y)&=y\brac{\frac{\alpha}{2}+(1-\alpha)\brac{\frac{n}{n-1}}^2 (1-y)^2}.
\end{align}
Since the function $h_n=g_n/f_n$ is continuously differentiable on $(0,1)$ and converges uniformly to a continuous limit $h$ as $n\to \infty$, it follows easily that $\sup_{n,y\in[r_\alpha,1/2]}\abs{h_n'(y)} < \infty$ and for any $\epsilon >0$ we have $\inf_{n,z\in[r_\alpha,1/2-\epsilon]} h_n(z) >1$.
Furthermore, from the second statement of Lemma~\ref{lem:f_complete}, it follows that, for $\alpha < 1/3$, the transition rates of $A$ satisfy~\eqref{eq:cond_neg_drft} with $r_1=r_\alpha$ and $r_2=1/2$. Finally, condition~\eqref{eq:cond_lazyness} is satisfied for all sufficiently large $n$ since we have
\begin{align}
    q_+(a)+q_-(a)=n\frac{\alpha}{2}+n(1-\alpha)\brac{\frac{n}{n-1}}^2\brac{\frac{a}{n}}\brac{1-\frac{a}{n}}\leq \frac{n}{2},   
\end{align}
which follows since $y(1-y)\leq 1/4$ for all $y\in (0,1)$ and $n/(n-1)\leq \sqrt 2$ for all sufficiently large $n$. 
Thus, by Lemma~\ref{lem:exp_bound_gen} we conclude that,  for each $p\in (r_\alpha,1/2)$ we have $T_{1/2}=\Omega(\exp(\Theta( n)))$ with high probability whenever $A(0)=\ceil{np}$. 
For $p\in (1/2, 1-r_\alpha)$ the proof follows similarly by applying Lemma~\ref{lem:exp_bound_gen} to the chain $n-A$.
For values of $p$ in the set $[0,r_\alpha]\cup[1-r_\alpha,1]$, the same bound on $T_{1/2}$ holds, since to reach $\ceil{n/2}$ or $\floor{n/2}$ starting from state $\ceil{np}$ requires either passing through the states $\ceil{np'}$ where $p'\in (r_\alpha,1/2)$, or passing through the states $\ceil{np'}$ where $p'\in (1/2, 1-r_\alpha)$.
Finally,
by Corollary \ref{cor: TV distance bound}, for any $p\in [0,1/2)$ we have 
\begin{equation}
d(t)\geq \frac{1}{2}-\PP_{\ceil{np}}(T_{\frac{1}{2}}\leq t).
\end{equation}
Now, choosing $t=\Omega(\exp(\Theta(n)))$ in the above inequality gives
\begin{equation}
d(\Omega(\exp(\Theta(n))))\geq \frac{1}{2}-\PP_{\ceil{np}}(T_{\frac{1}{2}}\leq \Omega(\exp(\Theta(n))))\geq \frac{1}{2}-o(1)> \frac{1}{4}.
\end{equation}
This implies $t_{\text{mix}}= \Omega(\exp(\Theta(n)))$.

To prove the second statement of the theorem, we use Lemma~\ref{lem:log_bound_gen}. 
First we define the function $\Phi$ on the state-space of $A$ as
\begin{equation}
    \Phi(a)=\sbrac{a-\ceil{\frac{n}{2}}}_+\wedge \sbrac{a-\floor{\frac{n}{2}}}_-.
\end{equation}

Then, for all $a\in \{0,1,\ldots,n\}$ we have 
\begin{align}
\mc G_A \Phi(a)&\leq(q_+(a)-q_-(a))\indic{a > \ceil{\frac{n}{2}}}+(q_-(a)-q_+(a))\indic{a < \floor{\frac{n}{2}}}\\
&\overset{(a)}{=}n\brac{f_{n,\alpha}\brac{\frac{a}{n}}-f_{n,\alpha}\brac{\frac{1}{2}}}\indic{a > \ceil{\frac{n}{2}}}+
n\brac{f_{n,\alpha}\brac{\frac{1}{2}}-f_{n,\alpha}\brac{\frac{a}{n}}}\indic{a < \floor{\frac{n}{2}}}\\
&\overset{(b)}{\leq} -nc(n,\alpha)\brac{\brac{\frac{a}{n}-\frac{1}{2}}\indic{a > {\frac{n}{2}}}+\brac{\frac{1}{2}-\frac{a}{n}}\indic{a 
<{\frac{n}{2}}}}\\
&=-c(n,\alpha)\Phi(a),
\end{align}
where $(a)$ follows from the definition of $f_{n,\alpha}$ and the fact that $f_{n,\alpha}(1/2)=0$, and (b) follows from the last statement of Lemma~\ref{lem:f_complete}. Now, for $\alpha> 1/3$, we have $c(n,\alpha)=(3\alpha-1)/2-o(1)>0$ for all sufficiently large $n$.
Hence, we can apply Lemma~\ref{lem:log_bound_gen} to $A$ with the above choice of $\Phi$ and choosing $\mc{Y}_0=\{\ceil{n/2},\floor{n/2}\}$. Furthermore, since $\Phi_{\min}\geq 1$ and $\Phi(a)\leq  n$ for any $a\in\{0,1,\ldots,n\}$, we conclude that
\begin{equation}
    \EE_{a}[T_{1/2}]\leq \frac{1+\log\brac{\Phi(a)}}{c(n,\alpha)}\leq \frac{1+\log n}{c(n,\alpha)},
\end{equation}
for any $a\in \{0,1,\ldots,n\}$ which shows that $\EE_{a}[T_{1/2}]=O(\log n)$. To prove the bound on $t_{\text{mix}}$ we use the coupling lemma for Markov chain mixing (see, for example, Theorem~5.4 of~\cite{peres_book}), according to which the distance $d(t)$ can be bounded above as
\begin{equation}
d(t) \leq \max_{a,\bar a\in \{0,1,\ldots,n\}}\PP_{a,\bar a}[A(t) \neq \bar A(t)]=\max_{a,\bar a\in \{0,1,\ldots,n\}}\PP_{a,\bar a}[T_{\text{couple}}>t],
\label{eq:dt_bound}
\end{equation}
where $A$ and $\bar A$ are two independent birth-death chains on $\{0,1,\ldots,n\}$, each evolving according to the transition rates~\eqref{eq:up}-\eqref{eq:down}, and
$T_{\text{couple}}=\inf\{t\geq 0: A(t)=\bar A(t)\}$ is their coupling time. To obtain a bound on $T_{\text{couple}}$, we define the function $\Phi$ on the state-space
of the joint chain $(A,\bar A)$ as $\Phi(a,\bar a)=\abs{\bar a-a}$.
If $\mc G_{A,\bar A}$ denotes the generator of the joint chain $(A,\bar A)$, then it is easy to see
that
\begin{align}
\mc G_{A,\bar A}\Phi(a,\bar a)&\leq(q_+(\bar a)+q_-(a)-q_-(\bar a)-q_+(a))\indic{\bar a > a}\nonumber\\
&\hspace{2em}+(q_+(a)+q_-(\bar a)-q_-(a)-q_+(\bar a))\indic{a > \bar a}\\
&=n(f_{n,\alpha}(\bar a/n)-f_{n,\alpha}(a/n))\indic{\bar a > a}+n(f_{n,\alpha}(a/n)-f_{n,\alpha}(\bar a/n))\indic{a > \bar a}\\
&\leq -c(n,\alpha)\Phi(a,\bar a),
\end{align}
where the last line follows from the last statement of Lemma~\ref{lem:f_complete}.
Again, since $c(n,\alpha)>0$ for all sufficiently large $n$, applying Lemma~\ref{lem:log_bound_gen} with $\Phi$ chosen as above,  $\mc{ Y}_0=\{(a,\bar a)\in \{0,1,\ldots,n\}\times \{0,1,\ldots,n\}: a=\bar a\}$, and $\Phi_{\min}=1$, we obtain
\begin{align}
\EE_{a,\bar a}[T_{\text{couple}}]
\leq \frac{1+\log\abs{\bar a-a}}{c(n,\alpha)}\leq \frac{1+\log n}{c(n,\alpha)}.
\end{align}
Hence, using the Markov inequality it follows that
\begin{equation}
    \PP_{a,\bar a}[T_{\text{couple}} > t] \leq \frac{\EE_{a,\bar a}{[T_{\text{couple}}]}}{t}\leq \frac{1+\log n}{tc(n,\alpha)}.
\end{equation}
Therefore, using~\eqref{eq:dt_bound} we have 
\begin{equation}
d(t) \leq \frac{1+\log n}{tc(n,\alpha)},
\end{equation}
which implies that for $t = \frac{4(1+\log n)}{c(n,\alpha)}$, we have $d(t)\leq \frac{1}{4}$. This shows that $t_{\text{mix}}\leq \frac{4(1+\log n)}{c(n,\alpha)}$, which completes the proof of the second statement of Theorem~\ref{thm:complete} since $c(n,\alpha)=(3\alpha-1)/2+o(1)$.\qed

{\em Remark on the 2$k$-choices dynamics}: We also remark here that the results established for the 2-choices dynamics with node failures on complete graphs can be extended to the $2k$-choices model, where, instead of sampling only two neighbours, an updating node samples $2k$ neighbours and takes the majority of this group and themself. The reason for $2k$ (even) choices is so that there are no split majorities. The rest of the model remains the same. For this dynamics, the drift function $f_{2k,n,\alpha}$ (analogous to the drift function $f_{n,\alpha}$) can be analysed through its limit $f_{2k,\alpha}$, given by
\begin{equation}
    f_{2k, \alpha}(y) =  \frac{\alpha}{2}(1-2y) + (1-\alpha) \sum_{r=k+1}^{2k}\binom{2k}{r}\left[y^r (1-y)^{2k-r+1} - (1-y)^r y^{2k-r+1}\right].
\end{equation}
It has been shown (in Appendix~\ref{proof:2k}) that this function satisfies analogous properties as stated in Lemma~\ref{lem:f_complete} for the 2-choices model. The only real difference is that the threshold value where $f_{2k,\alpha}$ changes from having three distinct roots to just one root in $[0,1]$ is given by
\begin{equation}
    \alpha_{2k} = 1- \frac{2^{2k}}{(2k+1) \binom{2k}{k}}.
\end{equation}
Hence, as expected, the threshold value increases with $k$, making the network more robust to failures.

\section{Extension to General Graphs}
\label{Section: General graphs}

In this section, we analyse the 2-choices model with node failures on general graphs.
To state the main result of this section, we need to introduce some additional notations. 
The adjacency matrix of the underlying graph $G=(V,E)$ is denoted by $\mf M=(M_{v,u})_{v,u\in V}$ where $M_{v,u}=M_{u,v}=\indic{\{u,v\} \in E}$ for all $u,v \in V$. We define the diagonal matrix $\mf D\in \RR^{n\times n}$ with entries $D_{v,u} = d_v \indic{v=u}$ for all $u,v \in V$. The scaled adjacency matrix $\mf{\bar M}$ is defined as $\mf{\bar M}=\mf{D}^{-1}\mf{M}$. 
We define the inner product $\langle \cdot, \cdot\rangle_{\mf D}$ between any two real vectors $\mf x, \mf y \in \RR^n$ as $\langle \mf x, \mf y \rangle_{\mf D}=\mf x^\top \mf D \mf y$ where $\cdot^\top$ denotes transpose operation. 
It is easy to verify that $\mf{\bar M}$ defines a self-adjoint linear operator on the Hilbert space $(\RR^n,\langle \cdot, \cdot\rangle_{\mf D})$, i.e., for any $\mf x, \mf y \in \RR^n$ we have $\langle \mf x, \mf{\bar M}\mf y\rangle_{\mf D}=\langle \mf{\bar M}\mf x, \mf y\rangle_{\mf D}$. Hence, by the spectral theorem, the eigenvalues of $\mf{\bar M}$ are real and the normalised eigenvectors of $\mf{\bar M}$ form an orthonormal basis for $\RR^n$. We denote the eigenvalues of $\mf{\bar M}$ by $\lambda_1,\lambda_2,\ldots,\lambda_n$ with the ordering $\lambda_1\geq \lambda_2\geq\ldots\geq \lambda_n$ and denote the normalised eigenvector corresponding to $\lambda_i$ by  $\mf f_i$ for each $i=1,2,\ldots,n$. It is easy to verify that $\lambda_1=1$ and 
\begin{equation}
    \mf f_1=\frac{1}{\sqrt{\vol(V)}}\mf 1,
    \label{eq:f1}
\end{equation}
where for any subset $S\subseteq V$, $\vol(S)$ denotes the sum of the degrees of the nodes in $S$, i.e., $\vol(S)=\sum_{v \in S} d_v$. We define $\lambda$ to be the second highest (in absolute value) eigenvalue of $\mf{\bar M}$, i.e., $\lambda=\max_{j \in \{2,3,\ldots,n\}}\abs{\lambda_j}<1$. The quantity $1-\lambda$ is often referred to as the {\em spectral gap} of the graph and indicates how well-connected the graph is (a higher spectral gap implies more connectivity). 

In this paper, we are interested in {\em expander graphs}, defined as a family of graphs with increasing $n$ such that the spectral gap remains bounded away from $0$, i.e, $\lambda < 1$ for all $n$. This definition includes both {\em dense graphs} where $\lambda=o(1)$ and {\em sparse graphs} where $\lambda=O(1)$. Examples of dense expander graphs include complete graphs where $\lambda=1/(n-1)$, \ErdosRenyi random graphs with an average degree $\Omega(\log n)$ where $\lambda=O(1/\sqrt{n})$ with high probability, and random $d$-regular graphs with $d=\omega(1)$ where $\lambda=O(1/\sqrt{d})$ with high probability. Examples of sparse expander graphs include random $d$-regular graphs with $d=O(1)$ where $\lambda=(2\sqrt{d-1}+o(1))/d$ with high probability. For our results to hold, we require the graph to be a member of such an expander family of graphs and to have sufficiently homogeneous degree distribution across nodes. To be more specific, we require that $\lambda \in (0,\bar \lambda)$ and $d_{\max}/d_{\min} < 1+\epsilon$ for some small $\bar \lambda\in (0,1)$ and small $\epsilon>0$ where $d_{\max}$ and $d_{\min}$, respectively, denote the maximum and minimum degrees of nodes in $G$. However, to state our result, instead of using $\bar \lambda$ and $\epsilon$, it would be more convenient to use an upper bound $L_{\max}$ on $(1+\lambda)^2\brac{{d_{\max}}/{d_{\min}}}^3$, and a lower bound $L_{\min}$ on $(1-\lambda)^2\brac{{d_{\min}}/{d_{\max}}}^3$. Note that for a graph with a large spectral gap and sufficiently homogeneous degrees, both $L_{\max}$ and $L_{\min}$ will be close to $1$.

The main result, which is an analogue of Theorem~\ref{thm:complete}, is given below.

\begin{theorem}
\label{thm:gen_graphs}
 
Consider an expander family of graphs satisfying $(1+\lambda)^2\brac{{d_{\max}}/{d_{\min}}}^3\leq L_{\max}$ and $(1-\lambda)^2\brac{{d_{\min}}/{d_{\max}}}^3\geq L_{\min}$, where $L_{\max}$ and $L_{\min}$ can possibly depend on $n$ and satisfy $L_{\max}> L_{\min} >0$. Furthermore, define $\Sigma_L=L_{\max}+L_{\min}$ and $\Delta_L=L_{\max}-L_{\min}$.
Then, the following results hold.
\begin{enumerate}
    \item  Let $0 < \alpha < \frac{\Sigma_L(1-K_L^{1/3})}{4+\Sigma_L(1-K_L^{1/3})}$. Then, for any graph with  $K_L:=\frac{27\Delta_L^2}{4\Sigma_L^2} < 1$, we have $t_{\text{mix}}= \Omega(\exp(\Theta(n)))$ and there exists $r\equiv r(\alpha,\Sigma_L,\Delta_L)\in(0,1/2)$ such that (i) $r\uparrow 1/2$ as $\Delta_L \to 0$, and (ii) with high probability $T_{1/2}= \Omega(\exp(\Theta(n)))$ given that $A(\mf X(0))=\ceil{np}$ for some $p\in [0,r)\cup(1-r,1]$.
    
    \item If $c(\alpha,L):={\alpha(4+\Sigma_L)-\Sigma_L} >0$, then the hitting time $\tau$ to the set 
    \begin{align*}
        \cbrac{\mf x: A(\mf x)=\ceil{\frac{n}{2} \brac{1+\epsilon_L}} \text{ or } A(\mf x)=\floor{\frac{n}{2} \brac{1-\epsilon_L}}},
    \end{align*}
    where $\epsilon_L=(1-\alpha)\frac{\Delta_L}{c(\alpha,L)}$ satisfies $\EE_{\mf x}[\tau]=O(\log n)$ for any initial state $\mf x \in \mc{X}$.
\end{enumerate}
\end{theorem}

{\em Discussion on the main result:} The above theorem is the analogue of Theorem~\ref{thm:complete} for general graphs. Note that the condition $K_L < 1$ in the first statement of the theorem requires $\Delta_L$ to be sufficiently small, which only occurs when $\lambda$ is small and $d_{\max}/d_{\min}$ is close to $1$.
Hence, the first statement of the theorem implies that as long as $\lambda$ is sufficiently small and $d_{\max}/d_{\min}$ is sufficiently close to $1$, there exists a threshold value for the failure probability $\alpha$, below which, the noisy 2-choices process takes an exponentially long time to mix and is therefore robust to node failures. The second statement characterises a set of values for $\alpha$ such that, for any given graph,
the chain reaches the region $\left(\frac{n}{2}(1-\epsilon_L),\frac{n}{2}(1+\epsilon_L)\right)$ in $O(\log n)$ time on average. Note that, unlike complete graphs, the second statement does not characterise the exact hitting time to the $n/2$ state. This is because providing a bound on the exact hitting time to $n/2$ would require making additional assumptions on the graph structure to obtain more refined bounds on the transition rates at states near $n/2$. However, even without any assumption on the graph structure, we are able to establish that when $\Delta_L$ is sufficiently small, $\epsilon_L$ is also small, implying that, for any given graph, the chain becomes sufficiently close to $n/2$ in $O(\log n)$ time, as long as $\alpha$ is above a threshold determined only by $\lambda$ and the ratio $d_{\max}/d_{\min}$ for that graph. 
It is of particular interest to derive the phase transition threshold for graphs that are asymptotically regular and dense, i.e., the properties $\lambda=o(1)$ and $d_{\max}/d_{\min}=1+o(1)$ are satisfied. For such graphs, we can choose $L_{\max}=1+o(1)$ and $L_{\min}=1-o(1)$ in the above theorem, which results in $\Sigma_L=2+o(1)$ and $\Delta_L=o(1)$. Plugging the above values of $\Sigma_L$ and $\Delta_L$ into Theorem~\ref{thm:gen_graphs}, we immediately obtain the following corollary.

\begin{corollary}
\label{cor:dense_graphs}
For an expander family of graphs satisfying $\lambda=o(1)$ and $d_{\max}/d_{\min}=1+o(1)$, the following results hold. 
\begin{enumerate}
    \item Fix $0< \alpha < \frac{1}{3}$ and choose any $\epsilon >0$. Then, $t_{\text{mix}}= \Omega(\exp(\Theta(n)))$ and for any $p\in [0,1/2-\epsilon)\cup(1/2+\epsilon,1]$ we have $T_{1/2}= \Omega(\exp(\Theta(n)))$ with high probability given that $A(\mf X(0))=\ceil{np}$.
    
    \item Fix $\alpha>1/3$ and choose any $\epsilon >0$. Then, the hitting time $\tau$ to the set $$\cbrac{\mf x: A(\mf x)=\ceil{\frac{n}{2} \brac{1+\epsilon}} \text{ or } A(\mf x)=\floor{\frac{n}{2} \brac{1-\epsilon}}},$$ satisfies $\EE_{\mf x}[\tau]=O(\log(n))$ for any initial state $\mf x \in \mc{X}$.
\end{enumerate}
\end{corollary}

\begin{figure}
    \centering
    \includegraphics[width=0.8\linewidth]{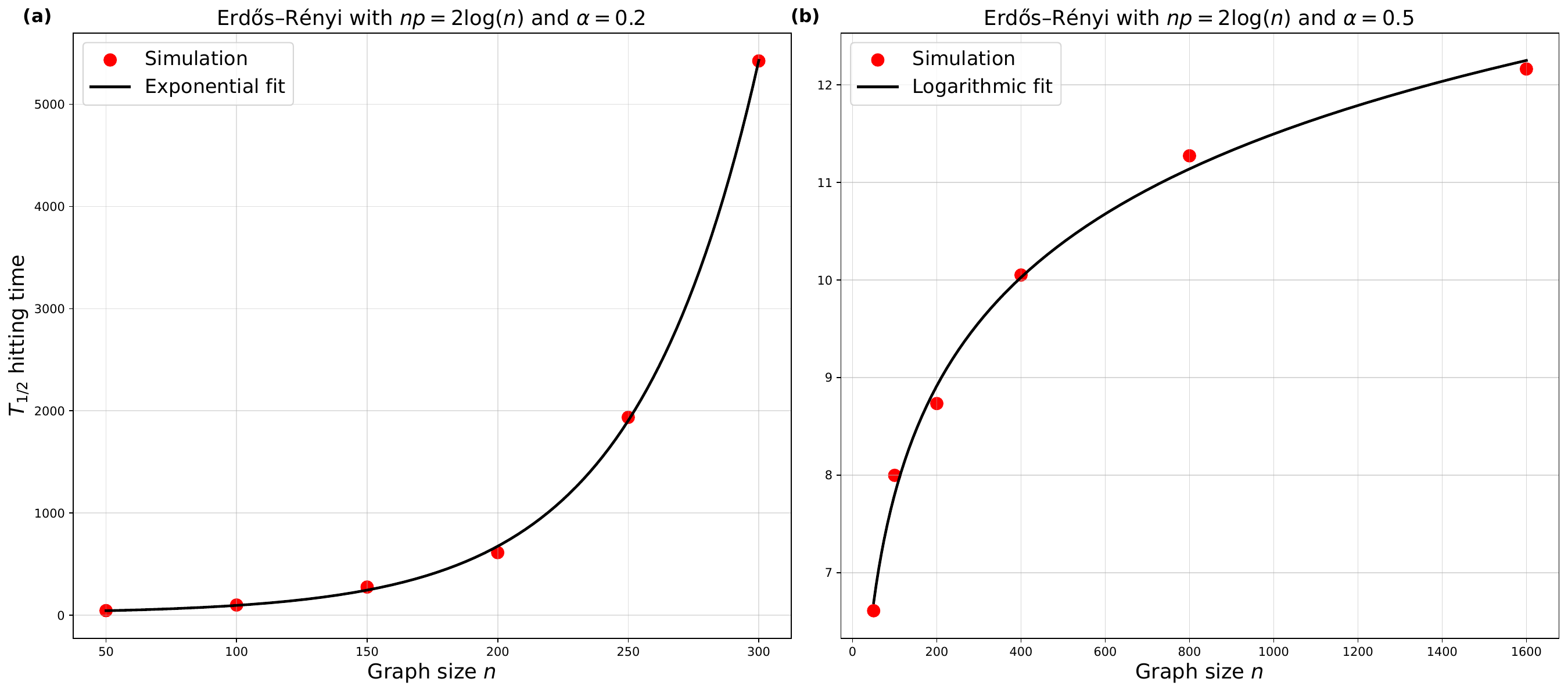}
    \caption{Plots showing $T_{1/2}$ as a function of network size for \ErdosRenyi random graphs with mean degree $2\log(n)$ for failure probabilities both below and above the critical threshold of $1/3$. A number of runs were performed on a fixed instance of the graph starting at state $\mf 0$, and the results were averaged. This was then repeated for a number of different \ErdosRenyi graph realisations with the same parameters and then the results were further averaged across these realisations. \textbf{(a)} $\alpha=0.2$. The number of simulations per graph was 10 and 5 graphs were generated per $n$ value. \textbf{(b)} $\alpha=0.5$. The number of simulations per graph was 30 and 10 graphs were generated per $n$ value.}
    \label{fig:ER graphs}
\end{figure}

\begin{figure}
    \centering
    \includegraphics[width=0.8\linewidth]{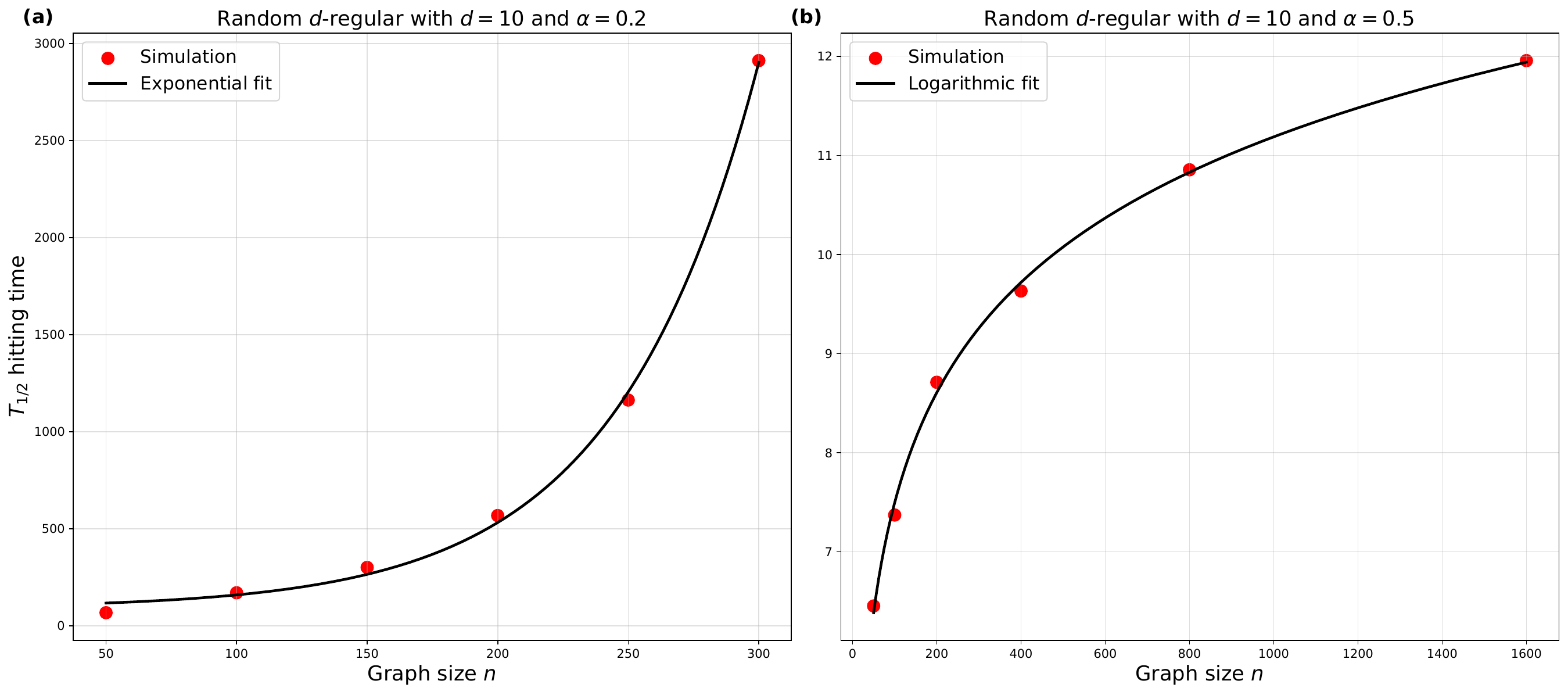}
    \caption{Plots showing $T_{1/2}$ as a function of network size for random $d$-regular graphs with $d=10$ for different values of the failure probability $\alpha$. A number of runs were performed on a fixed instance of the graph starting at state $\mf 0$, and the results were averaged. This was then repeated for a number of different random $d$-regular graph realisations with the same parameters and then the results were further averaged across these realisations. \textbf{(a)} $\alpha=0.2$. The number of runs per graph was 20 and 5 graph instances were generated per $n$ value. \textbf{(b)} $\alpha=0.5$. The number of runs per graph was 50 and 4 graphs were generated per $n$ value.}
    \label{fig:Dreg graphs}
\end{figure}

Hence, for asymptotically regular, dense graphs, the 2-choices dynamics essentially behaves in the same way as in complete graphs. For $\alpha <1/3$, both the mixing time and hitting time are exponentially large in $n$. In contrast, for $\alpha >1/3$ the chain reaches very close to the $n/2$ state in only logarithmic time. As mentioned earlier,  many important classes of graphs satisfy the properties of being dense and asymptotically regular. This includes \ErdosRenyi graphs with an average degree strictly above $\log n$, which gives asymptotically regular, dense graphs that are also strongly connected almost surely. In Figure~\ref{fig:ER graphs}, we plot the hitting time  to the state $n/2$ for \ErdosRenyi graphs as a function of $n$ for different values of $\alpha$. We clearly see that the results match our theory.

It is important to note that our results also hold for sparse graphs where $\lambda=O(1)$ and $d_{\max}/d_{\min}-1$ is small but constant. These are graphs where the average degree of a node remains bounded (does not scale with $n$). For example, for random $d$-regular graphs with $d=200$, the value of $\lambda$ can be computed using $\lambda=(2\sqrt{d-1})/d$, giving $\lambda\approx 0.14$. This allows us to calculate the threshold for metastability as $\alpha=0.09$. This shows that, for all $\alpha< 0.09$, the noisy 2-choices dynamics on a random regular graph with $d=200$ is robust to node failures. However, we point out that, unlike for dense graphs (where our bounds are asymptotically tight), our bounds for sparse graphs are pessimistic in that the 2-choices dynamics could be robust even for $\alpha$ values above the threshold provided in Theorem~\ref{thm:gen_graphs}. In Figure~\ref{fig:Dreg graphs}, we plot $T_{1/2}$ as a function of $n$ for random $d$-regular graphs. We see a phase transition similar to that in complete graphs.

\textit{Sketch proof of Theorem~\ref{thm:gen_graphs}:} The proof of Theorem~\ref{thm:gen_graphs} is similar in flavour to the proof of Theorem \ref{thm:complete} for complete graphs. The main difference though is that for the process $A(t)\equiv A(\mf X(t))$, tracking the number of nodes with opinion $1$ in the network is no longer Markovian. Indeed, the transition rates of $A(\mf X(t))$ at a given state $\mf X(t)=\mf x \in \mc{X}$ for a general graph are  given by
\begin{align}
    q_+(\mf x)&=(n-A(\mf x))\frac{\alpha}{2}+(1-\alpha)\sum_{v\in B(\mf x)}\brac{\frac{d_v^A(\mf x)}{d_v}}^2,\label{eq:up_gen} \text{ and}\\
    q_-(\mf x)&=A(\mf x)\frac{\alpha}{2}+(1-\alpha)\sum_{v\in A(\mf x)}\brac{\frac{d_v^B(\mf x)}{d_v}}^2,\label{eq:down_gen}
\end{align}
respectively. The above demonstrates how the transition rates now depend on the detailed state $\mf x$, as opposed to only $A(\mf x)$.  Despite this difficulty, we show that it is still useful to track $A(\mf X(t))$, since we are able to sandwich the transition rates of $A(\mf X(t))$ between a `slower' and a `faster' birth-death Markov chain, for which our previous results (specifically, Lemma~\ref{lem:exp_bound_gen} and Lemma~\ref{lem:log_bound_gen}) in Section~\ref{Section: Complete graphs} are still applicable. The bounds for the transition rates that underpin the slower and faster chains are developed in Lemma~\ref{lem:spectral_bound} below, using the spectral properties of the graph.

\begin{lemma}
\label{lem:spectral_bound}
Let $G=(V,E)$ be a graph with spectral gap $1-\lambda$ and maximum and minimum degrees denoted by $d_{\max}$  and $d_{\min}$, respectively. Then, for any two disjoint subsets $S,T \subseteq V$ with $S\cup T = V$, we have
\begin{align}
      L_{\min}\frac{|S|^2 |T|}{|V|^2}& \leq  \sum_{v \in T}\brac{\frac{d_v^S}{d_v}}^2  \leq  L_{\max}\frac{|S|^2 |T|}{|V|^2},\label{eq:flux_BA}
\end{align}
where $L_{\max}=(1+\lambda)^2  \left(\frac{d_{\max}}{d_{\min}}\right)^3$ and $L_{\min}=(1-\lambda)^2 \left(\frac{d_{\min}}{d_{\max}}\right)^3$.
\end{lemma}

From the above lemma and the transition rates in~\eqref{eq:up_gen}-\eqref{eq:down_gen}, it follows immediately that for any given graph $G$ we have
\begin{align*}
    \ubar q_+(A(\mf x)) &\leq q_+(\mf x) \leq \bar q_+(A(\mf x))\\
    \ubar q_-(A(\mf x)) &\leq q_-(\mf x) \leq \bar q_-(A(\mf x)),
\end{align*}
where 
\begin{align}
        & \bar q_+(a) = \frac{1}{2} \alpha (n - a)
       + (1-\alpha) L_{\max}\frac{a^2\bigl(n - a\bigr)}{n^2} \label{eq:bar_plus}\\
       & \ubar q_+(a) = \frac{1}{2} \alpha (n - a)
       + (1-\alpha) L_{\min}\frac{a^2\bigl(n - a\bigr)}{n^2} \label{eq:ubar_plus}\\
       & \bar q_-(a) = \frac{1}{2} \alpha a
       + (1-\alpha) L_{\max}\frac{a\bigl(n - a\bigr)^2}{n^2}\label{eq:bar_minus}\\
       & \ubar q_-(a) = \frac{1}{2} \alpha a
       + (1-\alpha) L_{\min}\frac{a\bigl(n - a\bigr)^2}{n^2}\label{eq:ubar_minus},
\end{align}
for all $a\in \{0,1,\ldots,n\}$. 

In order to prove our main results for general graphs, we use the simpler bounds above which depend on $\mf x$ only through $A(\mf x)$. This considerably simplifies our analysis, allowing us to focus on the dynamics defined by the bounds rather than the actual dynamics.
For example, in order to prove the exponential lower bound in the first part of Theorem~\ref{thm:gen_graphs}, we construct a birth-death process $\bar A$ on $\{0,1,\ldots,n\}$ with birth and death rates given by the bounds $\bar q_+(a)$ and $\ubar q_-(a)$, respectively, in state $a\in \{0,1,\ldots,n\}$. In the lemma below, proved in Appendix~\ref{proof:coupling_bound}, we construct a coupling between the original process $\mf X$ and the  process $\bar A$ defined above such that, starting from the same initial number of nodes with opinion $1$, the process $\bar A$ always remains ahead of the process $\mf X$, i.e., $\bar A(t)\geq A(\mf X(t))$ for all $t\geq 0$.

\begin{lemma}
    \label{lem:coupling_bound}
    Let $\bar A=(\bar A(t) : t \geq 0)$ be a birth-death process on $\{0,1,\ldots,n\}$ with birth and death rates in state $a$ given by $\bar q_+(\bar a)$ and $\ubar q_-(\bar a)$, respectively, where $\bar q_+$ and $\ubar q_-$ are as defined in~\eqref{eq:bar_plus}-\eqref{eq:ubar_minus}. Define $\bar T_r = \inf\{ t\geq 0 : \bar A(t) = \ceil{nr} or \bar A(t)=\floor{nr}\}$. Suppose $\bar A(0) = A(\mf X (0)) = a$, where $a \in \{0,1,\ldots,n\}$. Then, there exists a coupling between $\bar A$ and $\mf X$ such that  $\bar A(t) \geq A(\mf X(t))$ for all $t \geq 0$. Hence, for any $r\geq a/n$, we have $T_{r} \geq_{st}\bar T_{r}$.
\end{lemma}

Therefore, according to the above lemma, to establish a lower bound on $T_{1/2}$, it is sufficient to establish the same lower bound for $\bar T_{1/2}$. To obtain such bounds on the hitting and mixing times for the processes defined by the transition rates in~\eqref{eq:bar_plus}-\eqref{eq:ubar_minus}, we need the properties of the drifts defined by these functions. The next lemma, which is analogous to Lemma~\ref{lem:f_complete}, precisely does this.

\begin{lemma}
\label{lem:f_gen}
    Define $F_{\alpha,L}:\RR\to \RR$ as $F_{\alpha, L}(y) = (1-2y)\brac{y^2-y+\frac{\alpha}{\Sigma_L(1-\alpha)}}+\frac{\Delta_L}{4\Sigma_L}.$
   Then, $F_{\alpha,L}$ has the following properties.
    \begin{enumerate}
    \item For each $a\in \{0,1,\ldots,n\}$, it holds that 
    \begin{equation}
        \bar q_+(a)-\ubar q_-(a)=\bar q_-(n-a)-\ubar q_+(n-a)\leq n(1-\alpha)\frac{\Sigma_L}{2}F_{\alpha,L}\brac{\frac{a}{n}},
    \end{equation}
    where $\bar q_+, \ubar q_-, \bar q_-, \text{ and }\ubar q_+$ are as defined in~\eqref{eq:bar_plus}-\eqref{eq:ubar_minus}.
    \item Assume $K_L=\frac{27\Delta_L^2}{4 \Sigma_L^2} < 1$ and $\alpha < \frac{\Sigma_L(1-K_L^{1/3})}{4+\Sigma_L(1-K_L^{1/3})}$. Then, there are exactly two roots $r_1,r_2$ (with $r_1 < r_2$) of $F_{\alpha,L}$ in the interval $(0,1/2)$ and $F_{\alpha,L}(y) <0$ for all $y\in (r_1,r_2)$.
    
    \item If $c(\alpha,L):=(4+\Sigma_L)\alpha -\Sigma_L >0 $, then for all $y\in (1/2,1)$ we have 
    \begin{equation}
        F_{\alpha,L}(y)\leq \frac{2c(\alpha,L)}{\Sigma_L(1-\alpha)}\brac{\frac{1}{2}\brac{1+\frac{(1-\alpha)\Delta_L}{c(\alpha,L)}}-y}.
        \label{eq:gen_drift_bound}
    \end{equation}
\end{enumerate}
\end{lemma}

We now have all the ingredients required to prove the main result of this section. We proceed with the detailed proof below.

{\em Proof of Theorem~\ref{thm:gen_graphs}}: For $p\in [0,1/2)$, we establish the lower bound on $T_{1/2}$ in the first statement of the theorem  by applying  Lemma~\ref{lem:exp_bound_gen} to the process $\bar A$ defined in
Lemma~\ref{lem:coupling_bound}. For $p\in (1/2,1]$, the proof follows similarly by considering the process $n-\bar A$ instead of $\bar A$. Note that since $T_{1/2}\geq_{st}\bar T_{1/2}$ by Lemma~\ref{lem:coupling_bound}, it is sufficient to establish the lower bound for $\bar T_{1/2}$. From the first statement of Lemma~\ref{lem:f_gen}, it follows that, for the process $\bar A$, the drift satisfies
$\bar q_+(a) - \ubar q_-(a)\leq n(1-\alpha)\frac{\Sigma_L}{2}F_{\alpha,L}\brac{\frac{a}{n}},$
where $F_{\alpha,L}$ is as defined in Lemma~\ref{lem:f_gen}.
From the second statement of Lemma~\ref{lem:f_gen}, we conclude that the drift of $\bar A$ is negative in the interval $(r_1,r)$, where $r_1$ and $r$ (with $r_1 < r$) are the only two roots of $F_{\alpha,L}$ in the region $(0,1/2)$. From the definition of $F_{\alpha, L}$, it is clear that $r\to 1/2$ as $\Delta_L\to 0$. Furthermore, we have
\begin{align}
\bar q_+(a) + \ubar q_-(a)&=
          \frac{1}{2}\alpha n
           + (1-\alpha)\frac{a\bigl(n - a\bigr)}{n^2}(L_{\max}a+L_{\min}(n-a))\nonumber\\
&\leq \frac{1}{2}\alpha n
           + n(1-\alpha)L_{\max}\frac{a\bigl(n - a\bigr)}{n^2}\leq n\brac{\frac{1}{2}\alpha+(1-\alpha)\frac{L_{\max}}{4}}<  n,
\end{align}
where the last inequality follows since $\alpha < \frac{\Sigma_L}{4+\Sigma_L}< \frac{L_{\max}}{2+L_{\max}}< \frac{4-L_{\max}}{2-L_{\max}}$. Furthermore, it is easy to see that $h$ defined $h_n(z)=q_-(nz)/q_+(nz)$ is continuously differentiable and converges uniformly to some continuous function $h$ as $n\to \infty$. Hence, the process $\bar A$ satisfies the conditions of Lemma~\ref{lem:exp_bound_gen}, which imply the lower bound on $T_{1/2}$ stated in the first part of the theorem for $p \in (r_1,r)$. For any initial state $\ceil{np}$  with $p\in [0,r_1]$, the hitting time $\bar T_{1/2}$ is at least as large as the hitting time starting from states $\ceil{np}$ with $p \in (r_1,r)$. This shows that the hitting time bound holds for all $p\in [0,r)$. The proof of the mixing time bound now follows by applying Corollary~\ref{cor: TV distance bound} in a manner identical to that in the proof of Theorem~\ref{thm:complete}.

To prove the second part of the theorem, we apply Lemma~\ref{lem:log_bound_gen} to the process $\mf X$, with $\Phi: \mc{X} \to \RR$ defined as $\Phi(\mf x)=\sbrac{A(\mf x)-\ceil{\frac{n}{2}(1+\epsilon_L)}}_+\wedge\sbrac{A(\mf x)-\floor{\frac{n}{2}(1-\epsilon_L)}}_-$, where $\epsilon_L=(1-\alpha)\frac{\Delta_L}{c(\alpha,L)}$. Hence, $\tau$ is the first time $\Phi$ becomes zero.

For any $\mf x \in \mc{X}$ we have
\begin{align}
    G_{\mf X} \Phi(\mf x)&\leq(q_+(\mf x)-q_-(\mf x))\indic{A(\mf x) > {\frac{n}{2}}(1+\epsilon_L)}\nonumber\\
    &\hspace{2em}+(q_-(\mf x)-q_+(\mf x))\indic{A(\mf x) < {\frac{n}{2}}(1-\epsilon_L)}\nonumber\\
&\overset{(a)}{\leq}n(1-\alpha)\frac{\Sigma_L}{2}\left(F_{\alpha,L}\brac{\frac{A(\mf x)}{n}}\indic{A(\mf x) > {\frac{n}{2}}(1+\epsilon_L)}\right.\nonumber\\
&\hspace{2em}+\left.F_{\alpha,L}\brac{1-\frac{A(\mf x)}{n}}\indic{A(\mf x) < {\frac{n}{2}}(1-\epsilon_L)}\right)\nonumber\\
&\leq c(\alpha,L)\brac{\frac{n}{2}(1+\epsilon_L)-A(\mf x)}\indic{A(\mf x) > {\frac{n}{2}}(1+\epsilon_L)}\nonumber\\
&\hspace{2em}+c(\alpha,L)\brac{A(\mf x)-\frac{n}{2}(1-\epsilon_L)}\indic{A(\mf x) < {\frac{n}{2}}(1-\epsilon_L)}\nonumber\\
&=-c(\alpha,L)\Phi(\mf x),
\end{align}
where (a) follows from the first statement of Lemma~\ref{lem:f_gen}.
Hence, by Lemma~\ref{lem:log_bound_gen}, it follows that for any $\mf x$ we have
\begin{equation}
\EE_{\mf x}[\tau]\leq \frac{1+\log\Phi(\mf x)}{c(\alpha,L)}\leq \frac{1+\log n}{c(\alpha,L)},
\end{equation}
which completes the proof of the theorem.\qed


\section{Proofs of the general lemmas}
\label{sec:gen_proofs}

\subsection{Proof of Lemma \ref{lem:exp_bound_gen}}
\label{proof:exp_bound_gen}

By the hypothesis of the lemma, we have $h_n(z) > 1$ for all $z\in (r_1,r_2)$ and $h_n$ is a continuous function in $[r_1,r_2]$.

Next, note that the process $(\phi(Y(t)) : t \geq 0)$ where $\phi:\{0,1,\ldots,n\}\to \RR$ is defined as
\begin{equation}
k \mapsto \phi(k) = \sum_{t=1}^{k}\prod_{j=1}^{t-1}\frac{q_-
(j)}{q_+(j)},
\end{equation}
is a martingale with respect to the natural filtration generated by the process $Y$.
This follows by verifying that the generator $\mc{G}_Y$ of $Y$ applied to $\phi$ gives
$\mc G_{Y}\phi\equiv 0$. Hence, by the optional stopping theorem, we have $\EE_y[\phi(Y(\tau)] = \phi(y),$
for any (almost surely finite) stopping time $\tau$. Therefore, if we initialize the chain as $Y(0) = \ceil{np}$ for $p\in (r_1,r_2)$, and choose the stopping time $\tau$ as $\tau = T_{r_1} \wedge T_{r_2}$, which is almost surely finite by the ergodicity of $Y$ (which, in turn, follows from the condition $q_+(y)+q_-(y) > 0$ for all $y\in \{0,1,\ldots,n\}$), then
\begin{equation}
    \phi(\ceil{np}) = \phi(\ceil{nr_1})\PP_{\ceil{np}}(T_{r_1} < T_{r_2}) + \phi(\floor{nr_2})\PP_{\ceil{np}}(T_{r_1} > T_{r_2}).
\end{equation}
Now using $\PP_{\ceil{np}}(T_{r_1} < T_{r_2}) + \PP_{\ceil{np}}(T_{r_1} > T_{r_2}) = 1$, we obtain
\begin{equation}
    \PP_{\ceil{np}}(T_{r_2} < T_{r_1}) = \frac{\phi(\ceil{np}) - \phi(\ceil{nr_1})}{\phi(\floor{nr_2}) - \phi(\ceil{nr_1})}.
\end{equation}
Hence, using the definitions of $\phi$ and $h_n$, we obtain 
\begin{align}
\PP_{\ceil{np}}(T_{r_2} < T_{r_1}) 
&=\frac{\sum_{t=\ceil{nr_1}}^{\ceil{np}}\prod_{j=\ceil{nr_1}}^{t-1}h_n(j/n)}{\sum_{t=\ceil{nr_1}}^{\floor{nr_2}}\prod_{j=\ceil{nr_1}}^{t-1}h_n(j/n)}\nonumber\\
&\overset{(a)}{\leq} (\ceil{np}-\ceil{nr_1}+1)\frac{\prod_{j=\ceil{nr_1}}^{\ceil{np}-1}h_n(j/n)}{\prod_{j=\ceil{nr_1}}^{\floor{nr_2}-1}h_n(j/n)}\nonumber\\
&= \frac{\ceil{np}-\ceil{nr_1}+1}{\prod_{j=\ceil{np}}^{\floor{nr_2}-1}h_n(j/n)}\nonumber\\
&\leq\frac{{np}-{nr_1}+2}{\exp\brac{\sum_{j=\ceil{np}}^{\floor{nr_2}-1}\log(h_n(j/n))}}\nonumber\\
&\overset{(b)}{\leq}\frac{2n}{\exp\brac{n\int_{p}^{r_2}\log(h_n(z))dz+O(1)}}\\
& \overset{(c)}{\leq} \kappa_2 n \exp(-\kappa_1 n),
\label{eq: general hitting time bound}
\end{align}
where $\kappa_1$ and $\kappa_2$ are positive constants. In the above, (a) follows from the fact that $h(j/n) > 1$ for each $j \in (\ceil{nr_1}, \ceil{np})$, (b) follows from the facts that $n\geq 2$ and $(1/n)\sum_{j=\ceil{np}}^{\floor{nr_2}-1} \log h_n(j/n)= \int_{p}^{r_2} \log h_n(z)dz+O(1/n)$ for any function $h_n$ such that $\sup_{n,z\in [r_1,r_2]} \abs{h_n'(z)} < \infty$, and (c) follows from the fact that $\inf_{n,z\in[p,r_2-\epsilon]} h_n(z) > 1$ for any $\epsilon >0$.
Hence, if $Z_{r_1}$ denotes the number of returns to the state $\ceil{nr_1}$ before hitting $\floor{nr_2}$, then for any $i\geq1$ and all sufficiently large $n$ we have 
\begin{align}
    \PP_{\ceil{np}}(Z_{r_1} \geq i) &\geq (1 - \PP_{\ceil{np}}(T_{r_2} < T_{r_1}))^i\geq (1 - \kappa_2n\exp(-\kappa_1n))^i\geq 1 - i\kappa_2n\exp(-\kappa_1n),\label{eq:num_visit}
\end{align}
where the first inequality follows since the chain returns to state $\ceil{nr_1}$ every time $T_{r_1} < T_{r_2}$ starting from state $\ceil{pn}$, the second inequality follows from \eqref{eq: general hitting time bound}, and the final inequality follows since $(1-y)^i\geq 1-iy$ for $i \geq 1$ and $y\in (0,1)$. By choosing $i = \ceil{n \exp(\kappa_1 n/ 2)} \leq n \exp(\kappa_1 n / 2) + 1$ in~\eqref{eq:num_visit} we obtain
\begin{align}
    \PP_{\ceil{np}}[Z_{r_1} \geq \ceil{n \exp(\kappa_1 n/ 2)} ] \geq 1 - 2\kappa_2 n^2 \exp(-\kappa_1 n /2).
    \label{eq: number of visits to r_1 bound}
\end{align}
Note here that $T_{r_2} \geq \sum_{j=1}^{Z_{r_1}} M_{r_1, j}$, where $M_{r_1, j}$ denotes the time spent in the state $\ceil{nr_1}$ on the $j^{th}$ return to this state. Furthermore, each $M_{r_1, j}$ is an independent exponentially distributed random variable independent of $Z_{r_1}$, with a mean $1/q_+(\ceil{nr_1}) + q_-(\ceil{nr_1})\geq c/n$. We now apply a Chernoff bound to the sum of independent exponential random variables. For any non-negative integer valued function $\theta$, and all sufficiently large $n$, we have
\begin{align}
    \PP_{\ceil{np}}\left( \sum_{j=1}^{\theta(n)} M_{r_1, j} \leq \frac{\theta(n)}{n} \right)  \nonumber
    & = \PP_{\ceil{np}} \left( \exp \left( - \sum_{j=1}^{\theta(n)} M_{r_1, j}  \right)  \geq \exp \left( - \frac{\theta(n)}{n} \right) \right)\nonumber\\
    & \overset{(a)} \leq \exp \left(\frac{\theta(n)}{n} \right) \prod_{j=1}^{\theta(n)}\EE[\exp(-M_{r_1, j})]\nonumber\\
    & \overset{(b)} = \exp \left(\frac{\theta(n)}{n} \right) \left( \frac{q_+(\ceil{nr_1}) + q_-(\ceil{nr_1})}{1 + q_+(\ceil{nr_1}) + q_-(\ceil{nr_1})} \right)^{\theta(n)}\nonumber\\
    & \overset{(c)} \leq \exp \left(\frac{\theta(n)}{n} \right) \exp\left(-\frac{\theta(n)}{q_+(\ceil{nr_1}) + q_-(\ceil{nr_1})} \right)\nonumber\\
    & \overset{(d)} \leq \exp \left(\frac{\theta(n)}{n} \right) \exp\left(-\frac{c \theta(n)}{n} \right)\nonumber\\
    & = \exp \left(- (c-1) \frac{\theta(n)}{n} \right) = \exp \left(- \beta \frac{\theta(n)}{n} \right),
    \label{eq: Time spent bound}
\end{align}
where $\beta=c-1>0$. In the above, (a) follows from the Markov inequality and the independence of $\{M_{r_{1},j}\}_{j=1}^{\theta(n)}$, (b) follows from the moment generating function of an exponential random variable with mean $1/(q_+(\ceil{nr_1})+q_-(\ceil{nr_1}))$, (c) follows from the fact that $1+y\leq e^y$ for any $y\geq 0$, and (d) follows from ${1}/{(q_+(\ceil{nr_1}) + q_-(\ceil{nr_1}))}\geq c/n$. Hence, we have
\begin{align}
    \PP_{\ceil{np}}(T_{r_2} \geq \exp(n\kappa_1/2))\nonumber
    & \geq \PP_{\ceil{np}} \left( \sum_{j=1}^{Z_{r_1}} M_{r_1, j} \geq \exp(n\kappa_1/2) \right)\nonumber\\
    & \geq \PP_{\ceil{np}} \left( \sum_{j=1}^{Z_{r_1}} M_{r_1, j} \geq \exp(n\kappa_1/2), Z_{r_1} \geq \ceil{n \exp(n\kappa_1/2)}\right)\nonumber\\
    & \geq \PP_{\ceil{np}} \left( \sum_{j=1}^{\ceil{n \exp(n\kappa_1/2)}} M_{r_1, j} \geq \exp(n\kappa_1/2), Z_{r_1} \geq \ceil{n \exp(n\kappa_1/2)}\right)\nonumber\\
    & \geq\PP_{\ceil{np}} \left( \sum_{j=1}^{\ceil{n \exp(n\kappa_1/2)}} M_{r_1, j} \geq \frac{\ceil{n\exp(n\kappa_1/2)}}{n} \right) \nonumber\\
    &\hspace{10em}\times\PP_{\ceil{np}} \left( Z_{r_1} \geq \ceil{n \exp(n\kappa_1/2)}\right)\nonumber\\
    & \overset{(a)} \geq (1 - \exp (- \beta \exp(n \kappa_1 / 2)) (1 - 2\kappa_2 n^2 \exp(-\kappa_1 n /2)),
\end{align}
where (a) follows from \eqref{eq: number of visits to r_1 bound} and  \eqref{eq: Time spent bound}.
The statement of the lemma now easily follows from the RHS of the last inequality.\qed

\subsection{Proof of lemma \ref{lem:spectral_bound}}
\label{proof:spectral_bound}

We denote the indicator vectors of the two sets $S$ and $T$ by $\mf{1}(S)=(\indic{v\in S})_v$ and $\mf{1}(T)=(\indic{v\in T})_v$. Representing these in terms of the basis vectors $\mf f_1,\ldots,\mf f_n$, we have $\mathbf{1}(S) = \sum_{j=1}^n \beta_j(S)\,\mathbf{f}_j$ and $\mathbf{1}(T) = \sum_{j=1}^n \beta_j(T)\,\mathbf{f}_j,$
where for each $j \in \{1,2,\ldots,n\}$ and each $U \subseteq V$, we have $\beta_j(U)=\langle \mf{1}(U),\mf f_j\rangle_{\mf D}$. 

We first obtain bounds on the sum $\sum_{v \in T}\frac{\brac{d_v^S}^2}{d_v}$. Note that%
\begin{align}
\label{eq:final_lb}
\sum_{v \in T}\frac{\brac{d_v^S}^2}{d_v} =\vol(T) \sum_{v \in T}\frac{d_v}{\vol(T)}\brac{\frac{d_v^S}{d_v}}^2
&\overset{(a)}{\geq} \vol(T)\brac{\sum_{v \in T}\frac{d_v}{\vol(T)}\frac{d_v^S}{d_v}}^2\\
&= \frac{1}{\vol(T)}\brac{E(S, T)}^2\\
&\overset{(b)}{\geq} \frac{1}{\vol(T)}\brac{(1-\lambda)\frac{\vol(S)\vol(T)}{\vol(V)}}^2\\
&=(1-\lambda)^2\frac{\vol^2(S)\vol(T)}{\vol(V)},
\end{align}
where (a) follows by Jensen's inequality and (b) follows from Lemma~\ref{lem:expander_mixing}.
To derive an upper bound on $\sum_{v \in T}\frac{\brac{d_v^S}^2}{d_v}$, we note that
\begin{align}
\sum_{v \in V}\frac{\brac{d_v^S}^2}{d_v}= \norm{\mf{\bar M}\mf{1}(S)}_{\mf D}^2
=\sum_{j=1}^n \beta_j^2(S)\lambda_j^2 &= \beta_1^2(S) + \sum_{j=2}^n \beta^2_{j}(S)\lambda^2_j\\ 
&\leq \frac{\vol^2(S)}{\vol(V)}+\lambda^2\sum_{j=2}^n \beta_j^2(S)\\
&\leq \frac{\vol^2(S)}{\vol(V)}+\lambda^2\frac{\vol(S)\vol(T)}{\vol(V)},
\end{align}
where in the last two lines we have used Lemma~\ref{lem:indicator}.
Furthermore, we have
\begin{align}
\sum_{v \in S} \frac{\brac{d_v^S}^2}{d_v}= \vol(S)  \sum_{v \in S} \frac{d_v}{\vol(S)}\brac{\frac{d_v^S}{d_v}}^2
&\geq \vol(S)  \brac{\sum_{v \in S} \frac{d_v}{\vol(S)}\frac{d_v^S}{d_v}}^2\\
&=\frac{1}{\vol(S)}\brac{\vol(S)-\sum_{v \in S} d_v^T}^2\\
&=\frac{1}{\vol(S)}\brac{\vol(S)-E(S,T)}^2\\
&\geq \frac{1}{\vol(S)}\brac{\vol(S)-(1+\lambda)\frac{\vol(S)\vol(T)}{\vol(V)}}^2\\
&=\vol(S)\brac{\frac{\vol(S)}{\vol(V)}-\lambda \frac{\vol(T)}{\vol(V)}}^2.
\end{align}
Hence, we have
\begin{align}
\sum_{v \in T} \frac{\brac{d_v^S}^2}{d_v}&=\sum_{v \in V} \brac{\frac{d_v^S}{d_v}}^2-\sum_{v \in S} \brac{\frac{d_v^S}{d_v}}^2\nonumber\\
&\leq \frac{\vol^2(S)}{\vol(V)}+\lambda^2\frac{\vol(S)\vol(T)}{\vol(V)}- \vol(S)\brac{\frac{\vol(S)}{\vol(V)}-\lambda \frac{\vol(T)}{\vol(V)}}^2\nonumber\\
&=(1+\lambda)^2\frac{\vol^2(S)\vol(T)}{\vol^2(V)}.\label{eq:final_ub}
\end{align}
Using the fact that $d_{\min} \abs{U} \leq \vol(U)\leq d_{\max} \abs{U}$ for any $U\subseteq V$ and $\frac{1}{d_{\max}} \sum_{v \in T}\frac{\brac{d_v^{S}}^2}{d_v}\leq \sum_{v \in T}\brac{\frac{d_v^{S}}{d_v}}^2\leq \frac{1}{d_{\min}} \sum_{v \in T}\frac{\brac{d_v^{S}}^2}{d_v}$ in~\eqref{eq:final_ub} and~\eqref{eq:final_lb}, we can easily derive~\eqref{eq:flux_BA}. This completes the proof of the lemma.

\section{Conclusion}
\label{Section: Conclusion}

We have studied the 2-choices model under the effect of node failures. We assumed a stochastic model of failure in which a node fails with a constant probability $\alpha$ at each update instant.  Our results for complete graphs showed a sharp phase transition at $\alpha=1/3$, where the system undergoes a change from metastability  to rapid mixing.  
For general graphs, we showed that similar behaviour exists for both sparse and dense graphs, as long as they have sufficiently large spectral gaps and sufficiently homogeneous degree distributions. 

There are many open problems worth investigating. For example, our results do not apply to power-law graphs as they often have very heterogeneous degree distributions. Furthermore, the bound on the failure probability we obtain for sparse graphs is not asymptotically tight. Finding asymptotically tight bounds for sparse graphs will likely require more a careful analysis using local weak convergence results on graphs. Finally, our stochastic model of failure  does not capture scenarios where failure occurs to $\Theta(n)$ adversarially chosen  nodes per unit time.

\bibliographystyle{unsrt}
\bibliography{opinion}

@article{arpan_JSP_20,
  title={Voter and majority dynamics with biased and stubborn agents},
  author={Mukhopadhyay, Arpan and Mazumdar, Ravi R and Roy, Rahul},
  journal={Journal of Statistical Physics},
  volume={181},
  number={4},
  pages={1239--1265},
  year={2020},
  publisher={Springer}
}

@inproceedings{cooper_two_choices,
  title={The power of two choices in distributed voting},
  author={Cooper, Colin and Els{\"a}sser, Robert and Radzik, Tomasz},
  booktitle={International Colloquium on Automata, Languages, and Programming},
  pages={435--446},
  year={2014},
  organization={Springer}
}

@article{redner_two_choices,
  title={Dynamics of majority rule in two-state interacting spin systems},
  author={Krapivsky, Paul L and Redner, Sidney},
  journal={Physical Review Letters},
  volume={90},
  number={23},
  pages={238701},
  year={2003},
  publisher={APS}
}

@article{cruise2014probabilistic,
  title={Probabilistic consensus via polling and majority rules},
  author={Cruise, James and Ganesh, Ayalvadi},
  journal={Queueing Systems},
  volume={78},
  number={2},
  pages={99--120},
  year={2014},
  publisher={Springer}
}

@article{holley1975ergodic,
  title={Ergodic theorems for weakly interacting infinite systems and the voter model},
  author={Holley, Richard A and Liggett, Thomas M},
  journal={The Annals of Probability},
  pages={643--663},
  year={1975},
  publisher={JSTOR}
}

@article{clifford1973model,
  title={A model for spatial conflict},
  author={Clifford, Peter and Sudbury, Aidan},
  journal={Biometrika},
  volume={60},
  number={3},
  pages={581--588},
  year={1973},
  publisher={Oxford University Press}
}

@article{cox1989coalescing,
  title={Coalescing random walks and voter model consensus times on the torus in Z d},
  author={Cox, J Theodore},
  journal={The Annals of Probability},
  pages={1333--1366},
  year={1989},
  publisher={JSTOR}
}

@article{d2022phase,
  title={Phase transition of a nonlinear opinion dynamics with noisy interactions},
  author={d’Amore, Francesco and Clementi, Andrea and Natale, Emanuele},
  journal={Swarm Intelligence},
  volume={16},
  number={4},
  pages={261--304},
  year={2022},
  publisher={Springer}
}

@article{vieira2016phase,
  title={Phase transitions in the majority-vote model with two types of noises},
  author={Vieira, Allan R and Crokidakis, Nuno},
  journal={Physica A: Statistical Mechanics and its Applications},
  volume={450},
  pages={30--36},
  year={2016},
  publisher={Elsevier}
}

@article{cooper2013coalescing,
  title={Coalescing random walks and voting on connected graphs},
  author={Cooper, Colin and Elsasser, Robert and Ono, Hirotaka and Radzik, Tomasz},
  journal={SIAM Journal on Discrete Mathematics},
  volume={27},
  number={4},
  pages={1748--1758},
  year={2013},
  publisher={SIAM}
}

@inproceedings{nakata1999probabilistic,
  title={Probabilistic local majority voting for the agreement problem on finite graphs},
  author={Nakata, Toshio and Imahayashi, Hiroshi and Yamashita, Masafumi},
  booktitle={International Computing and Combinatorics Conference},
  pages={330--338},
  year={1999},
  organization={Springer}
}

@article{cruciani2021phase,
  title={Phase transition of the 2-Choices dynamics on core--periphery networks},
  author={Cruciani, Emilio and Natale, Emanuele and Nusser, Andr{\'e} and Scornavacca, Giacomo},
  journal={Distributed Computing},
  volume={34},
  number={3},
  pages={207--225},
  year={2021},
  publisher={Springer}
}

@article{chen2005majority,
  title={Majority rule dynamics in finite dimensions},
  author={Chen, P and Redner, S},
  journal={Physical Review E—Statistical, Nonlinear, and Soft Matter Physics},
  volume={71},
  number={3},
  pages={036101},
  year={2005},
  publisher={APS}
}

@inproceedings{doerr2011stabilizing,
  title={Stabilizing consensus with the power of two choices},
  author={Doerr, Benjamin and Goldberg, Leslie Ann and Minder, Lorenz and Sauerwald, Thomas and Scheideler, Christian},
  booktitle={Proceedings of the twenty-third annual ACM symposium on Parallelism in algorithms and architectures},
  pages={149--158},
  year={2011}
}

@inproceedings{vojnovic,
  title={Using three states for binary consensus on complete graphs},
  author={Perron, Etienne and Vasudevan, Dinkar and Vojnovic, Milan},
  booktitle={IEEE INFOCOM 2009},
  pages={2527--2535},
  year={2009},
  organization={IEEE}
}

@inproceedings{ghaffari2018nearly,
  title={Nearly-tight analysis for 2-choice and 3-majority consensus dynamics},
  author={Ghaffari, Mohsen and Lengler, Johannes},
  booktitle={Proceedings of the 2018 ACM Symposium on Principles of Distributed Computing},
  pages={305--313},
  year={2018}
}

@book{Kallenberg_book,
  author    = {Kallenberg, Olav},
  title     = {Foundations of Modern Probability},
  edition   = {2},
  publisher = {Springer-Verlag},
  address   = {New York},
  year      = {2002},
  series    = {Probability and Its Applications},
  isbn      = {0-387-95313-2}
}

@inproceedings{cooper2015fast,
  title={Fast consensus for voting on general expander graphs},
  author={Cooper, Colin and Els{\"a}sser, Robert and Radzik, Tomasz and Rivera, Nicolas and Shiraga, Takeharu},
  booktitle={International Symposium on Distributed Computing},
  pages={248--262},
  year={2015},
  organization={Springer}
}

@book{peres_book,
  title={Markov chains and mixing times},
  author={Levin, David A and Peres, Yuval},
  volume={107},
  year={2017},
  publisher={American Mathematical Soc.}
}

@article{noisy_voter_model,
  title={The noisy voter model},
  author={Granovsky, Boris L and Madras, Neal},
  journal={Stochastic Processes and their applications},
  volume={55},
  number={1},
  pages={23--43},
  year={1995},
  publisher={Elsevier}
}

@article{redner2019reality,
  title={Reality-inspired voter models: A mini-review},
  author={Redner, Sidney},
  journal={Comptes Rendus. Physique},
  volume={20},
  number={4},
  pages={275--292},
  year={2019}
}

@article{mukhopadhyay2024phase,
  title={Phase transitions in biased opinion dynamics with 2-choices rule},
  author={Mukhopadhyay, Arpan},
  journal={Probability in the Engineering and Informational Sciences},
  volume={38},
  number={2},
  pages={227--244},
  year={2024},
  publisher={Cambridge University Press}
}

@article{tay2008impact,
  title={On the impact of node failures and unreliable communications in dense sensor networks},
  author={Tay, Wee Peng and Tsitsiklis, John N and Win, Moe Z},
  journal={IEEE Transactions on Signal Processing},
  volume={56},
  number={6},
  pages={2535--2546},
  year={2008},
  publisher={IEEE}
}

@article{briat2023noise,
  title={Noise in biomolecular systems: Modeling, analysis, and control implications},
  author={Briat, Corentin and Khammash, Mustafa},
  journal={Annual Review of Control, Robotics, and Autonomous Systems},
  volume={6},
  number={1},
  pages={283--311},
  year={2023},
  publisher={Annual Reviews}
}

\newpage
\appendix
\section{Proof of Lemma~\ref{lem:stationary}}
\label{proof:stationary}

We construct a coupling between two Markov chains $\mf X$ and $\mf{\bar X}$, each evolving according to the transition rates given by~\eqref{eq:rates_graph}, but starting at complementary states, i.e., if $\mf X$ starts in state $\mf x$, then we start $\mf{\bar X}$ in state $\mf{\bar x}=\bm{1} - \mf x$. Under the coupling, for each node $v \in V$, the associated unit rate Poisson process triggering the updates of $v$ is the same in both $\mf X$ and $\mf{\bar X}$. Furthermore, if node $v$ performs an erroneous update in chain $\mf X$ choosing opinion $i\in \{0,1\}$, then the same update occurs in $\mf{\bar X}$. Similarly, if a node $v$ performs a two-choice update  by sampling neighbours $u$ and $u'$ in chain $\mf X$, then the same node $v$ also performs a two-choice update in chain $\mf{\bar X}$ by selecting the same two neighbours $u$ and $u'$. It is easy to see that while each chain individually obeys the transition rates given by~\eqref{eq:rates_graph}, their states remain complementary at all times, i.e., whenever $\mf X$ visits the state $\mf y$, $\mf{\bar X}$ visits the state $\mf{\bar y}=\bm{1} - \mf y$. Since both chains have the same transition rates, they have the same stationary distribution $\bs \pi$, and it therefore follows that $\bs \pi(\mf y)=\bs \pi(\bm{1} -\mf y)$ for all $\mf y \in \{0,1\}^n$. This proves the first part of the lemma. To prove the second part, we observe that 

\begin{align*}
\EE_{\mf X \sim \bs \pi}[A(\mf X)]&=\sum_{\mf x \in \{0,1\}^n} A(\mf x) \bs{\pi}(\mf x)\\
&=\sum_{\mf x:A(\mf x) < n/2} (A(\mf x) \bs{\pi}(\mf x)+A(\bm{1} -\mf x)\bs{\pi}(\bm{1} -\mf x))+ \frac{n}{2} \sum_{\mf x: A(\mf x)=n/2}\bs{\pi}(\mf x)\\
&=\sum_{\mf x:A(\mf x) < n/2} (A(\mf x)+A(\bm{1} -\mf x)) \bs{\pi}(\mf x)+ \frac{n}{2} \sum_{\mf x: A(\mf x)=n/2}\bs{\pi}(\mf x)\\
&=\sum_{\mf x:A(\mf x) < n/2} n \bs{\pi}(\mf x)+ \frac{n}{2} \sum_{\mf x: A(\mf x)=n/2}\bs{\pi}(\mf x)\\
&=\frac{n}{2}\sum_{\mf x:A(\mf x) \neq n/2} \bs{\pi}(\mf x)+ \frac{n}{2} \sum_{\mf x: A(\mf x)=n/2}\bs{\pi}(\mf x)=\frac{n}{2},
\end{align*}
where the equality in the third line follows from the first part of the lemma.

\section{Proof of Lemma~\ref{lem:f_complete}}
\label{proof:f_complete}

The first statement follows easily from the definition of $f_{n,\alpha}$.

To prove the second statement, note that $f_{n,\alpha}$ can be factorised as  follows
\begin{equation}
    f_{n,\alpha}(y)=(1-\alpha)\brac{\frac{n}{n-1}}^2(1-2y)\left(\frac{1}{2}\frac{\alpha}{1-\alpha}\brac{\frac{n-1}{n}}^2-y(1-y)\right).
    \label{eq:factored_f}
\end{equation}
The fourth factor in the above expression is a quadratic polynomial in $y$ with discriminant $1-\frac{2\alpha}{(1-\alpha)}\brac{1-\frac{1}{n}}^2$. 
Hence, this factor has two distinct real roots only when $\alpha < \frac{1}{1+2(1-\frac{1}{n})^2}$, which is implied by the condition $\alpha < 1/3$. Furthermore, since the sum of the roots is $1$ and their product is positive, both roots must lie in the range $(0,1)$. We denote the smaller of the two roots as $r_{n,\alpha}$, i.e. 
\begin{equation}
    r_{n,\alpha}=\frac{1}{2}\brac{1-\sqrt{1-\frac{2\alpha}{(1-\alpha)}\brac{1-\frac{1}{n}}^2}}.
\end{equation}
Hence, $f_{n,\alpha}(y) < 0$ for all $y \in (r_{n,\alpha},1/2)$.
The second part of the lemma now follows since $r_{n,\alpha}< r_{\alpha}$ for all $n$.

To prove the last statement,
we define $b_n = (n/(n-1))^2$ and rewrite $f_{n,\alpha}$ as
\begin{align}
    f_{n, \alpha}(y) = (1-2y)\frac{\alpha}{2}-(1-\alpha)b_ny(1-y)(1-2y).
\end{align}
The derivative $f_{n,\alpha}'$ of this function is given by 
\begin{align}
    f'_{n, \alpha}(y) = - [\alpha + (1-\alpha)b_n (6y^2-6y+1)].
\end{align}
The derivative above takes its maximum value in the interval $(0,1)$ at $y=1/2$ since
the second derivative given by
\begin{align}
    f''_{n, \alpha}(y) = -(1-\alpha)b_n (12y-6),
\end{align}
is zero when $y=1/2$, positive when $y < 1/2$, and negative when $y>1/2$. Hence, the maximum value of $f_{n,\alpha}'$ in $(0,1)$ is
\begin{align}
    f'_{n, \alpha}(1/2) =-\alpha + \frac{1}{2}(1-\alpha)b_n.
\end{align}
Hence, for $1\geq y \geq y'\geq 0$ we have 
\begin{align}
    f_{n, \alpha}(y) - f_{n, \alpha}(y') = \int_{y'}^y f'_{n, \alpha}(u)du \leq f_{n,\alpha}'(1/2) \int_{y'}^y du = -c(n,\alpha)(y-y'),
\end{align}
where $c(n,\alpha)=-f_{n,\alpha}'(1/2)=\frac{3\alpha-1}{2}-\frac{1}{2}(1-\alpha)(b_n-1)=\frac{3\alpha-1}{2}-o(1)$. This completes the proof of the last statement.\qed

\section{Proof of Lemma~\ref{lem:log_bound_gen}}
\label{proof:log_bound_gen}

To prove the lemma, we shall use Dynkin's formula~\cite{Kallenberg_book}, which is applicable to $Y$ since it is a CTMC with finite state-space, and therefore satisfies the Feller property. According to Dynkin's formula 
for any function $\Psi:\mc Y\to \RR$ and any stopping time $\tau'$ satisfying $\EE_{y}[\tau'] < \infty$ for all $y\in \mc{Y}$, we have 
\begin{equation}
    \EE_{y}[\Psi(Y(\tau'))]=\Psi(y)+\EE_{y}\sbrac{\int_{0}^{\tau'} \mc G_{Y} \Psi(Y(s))ds},
    \label{eq:dynkin}
\end{equation}
where $\mc G_Y$, as defined in the lemma's statement, is the generator of $Y$.
Let us define the function
$V:\RR_+\to \RR_+$ as
\begin{equation}
V(z)=\begin{cases}
     1+\log \brac{\frac{z}{\Phi_{\min}}}, &\text{ if } z > 0,\\
     0, &\text{otherwise.}
\end{cases}
\end{equation}
Now we apply Dynkin's formula to the function $\Psi=V \circ \Phi$ using the bounded stopping time $\tau_{z,j}=\tau_z\wedge j$, where $\Phi$ is defined in the statement of the lemma and $j\in \RR_+$. Before applying the formula, we note that $\Psi=V \circ \Phi \geq 0$, and for any $y\in B-\mc{ Y}_0$ we have
\begin{align}
\mc G_Y \Psi(y)&=\mc G_Y V(\Phi(y))\nonumber\\
&=\sum_{y'\neq y} q_{y,y'}(V(\Phi(y'))-V(\Phi(y)))\nonumber\\
&=\sum_{y'\neq y}q_{y,y'}\brac{\log \brac{\frac{\Phi(y')}{\Phi(y)}}\indic{y'\notin\mc{ Y}_0}-(1+\log \brac{\frac{\Phi(y)}{\Phi_{\min}}})\indic{y'\in \mc{ Y}_0}}\nonumber\\
&\leq\sum_{y'\neq y}q_{y,y'}\brac{\log \brac{1+\frac{\Phi(y')-\Phi(y)}{\Phi(y)}}\indic{y'\notin \mc{ Y}_0}-\indic{y'\in \mc{ Y}_0}}\nonumber\\
&\overset{(a)}{\leq}\sum_{y'\neq y}q_{y,y'}\brac{\frac{\Phi(y')-\Phi(y)}{\Phi(y)}\indic{y'\notin \mc{ Y}_0}+\frac{\Phi(y')-\Phi(y)}{\Phi(y)}\indic{y'\in \mc{ Y}_0}}\nonumber\\
&=\sum_{y'\neq y}q_{y,y'}\brac{\frac{\Phi(y')-\Phi(y)}{\Phi(y)}}\nonumber\\
&{=}\frac{\mc{G}_Y \Phi(y)}{\Phi(y)}\nonumber\\
&\overset{(b)}{\leq} -c\label{eq:drift_bound},
\end{align}
where $q_{y,y'}$ denotes the transition rate of $Y$ from state $y\in \mc Y$ to state $y'\in \mc{Y}$. In (a) we have used the fact that $\log (1+w) \leq w$ for $w> -1$ and $\Phi(y) \geq 0$ for all $y\in \mc{\bar Y}$, with equality holding iff $y\in \mc{ Y}_0$, and (b) follows since $\mc G_Y \Phi(y) \leq -c \Phi(y)$ by the lemma's hypothesis.
Hence, using~\eqref{eq:drift_bound} in Dynkin's formula we obtain
\begin{equation}
0\leq \EE_y[V(\Phi(Y(\sigma_j)))]\leq V(\Phi(y))-c\EE_y[\tau_{z,j}],
\end{equation}
which gives
\begin{equation}
\EE_y[\tau_{z,j}]\leq \frac{V(\Phi(y))}{c}=\frac{1+\log (\Phi(y)/\Phi_{\min})}{c}.
\end{equation}
Since the above is true for any $j$ and $\tau_{z,j}=\tau_z\wedge j \to \tau_z$ as $j \to \infty$, the result of the lemma follows by taking the limit of the above as $j\to \infty$ and using the dominated convergence theorem. 
\qed

\section{Expander Mixing Lemma}
\label{proof:expander_mixing}


The indicator vector of any subset of vertices can be written as a linear combination of the basis vectors $\mf f_1,\ldots,\mf f_n$ as follows
\begin{equation}
    \mf 1(U)=\sum_{j=1}^n\beta_j(U)\mf f_j,
\end{equation}
where $\beta_j(U)=\langle \mf{1}(U),f_j\rangle_{\mf D}$ for each $j \in \{1,2,\ldots,n\}$. Before stating the expander mixing lemma it will be helpful state some properties of the coefficients $\{\beta_j(U), j=1,2,\ldots,n\}$ for any $U\subseteq V$.

\begin{lemma}
\label{lem:indicator}
For any $U\subseteq V$, let $\mf 1(U)=\sum_{j=1}^n\beta_j(U)\mf f_j$. Then we have 
\begin{equation}
\beta_1(U)=\frac{\vol(U)}{\sqrt{\vol{V}}}.
\label{eq:beta_1}
\end{equation}
Furthermore, for any two disjoint subsets $S$ and $T$ of $V$ satisfying $S\cup T=V$ we have
$\beta_j(S)=-\beta_j(T)$ for all $j\in\{2,3,\ldots,n\}$. Furthermore, we have
\begin{equation}
    \sum_{j=2}^n \beta_j^2(S)=\sum_{j=2}^n \beta_j^2(T)=\frac{\vol(S)\vol(T)}{\vol(V)}.
\label{eq:beta_j}
\end{equation}
\end{lemma}

\begin{proof}
Since we have $$\mf f_1=\frac{1}{\sqrt{\vol(V)}}\mf 1,$$ and $\beta_1(U)=\langle \mf{1}(U),f_1\rangle_{\mf D}$,
\eqref{eq:beta_1} follows directly from the definition of the inner product $\langle\cdot,\cdot\rangle_{\mf D}$. Also note that for two disjoint sets $S$ and $T$ forming a partition of $V$ we have
\begin{align}
    \mf 1(S)+\mf 1(T)&=1, \quad \text{and}\label{eq:sum}\\
    \langle \mf 1(S),\mf 1(T)\rangle_D& =0\label{eq:ortho}.
\end{align}
Hence, from~\eqref{eq:sum} it follows that
\begin{align}
    \beta_j(S) &= \langle\mf{1}(S), \mathbf{f}_j \rangle_{\mf D}\\
    & = \langle \mf{1} - \mf{1}(T), \mathbf{f}_j\rangle_{\mf D}\\
    & = \langle \mf{1}, \mathbf{f}_j\rangle_{\mf D} - \langle\mf{1}(T), \mathbf{f}_j \rangle_{\mf D}\\
    &=\sqrt{\vol(V)}\langle \mathbf{f}_1, \mathbf{f}_j\rangle_{\mf D} - \beta_j(T),
\end{align}
where the last line follows from~\eqref{eq:f1}.
Hence, for $j \geq 2$, by the orthogonality of the eigenvectors we have $\beta_j(S)=-\beta_j(T)$.
Hence, from~\eqref{eq:ortho} we have
\begin{align}
    0=\sum_{j=1}^n \beta_j(S)\beta_j(T)=\beta_1(S)\beta_1(T)-\sum_{j=2}^n\beta_j^2(S)=\beta_1(S)\beta_1(T)-\sum_{j=2}^n\beta_j^2(T).
\end{align}
Hence, we have
\begin{equation}
\sum_{j=2}^n \beta_j^2(S)=\sum_{j=2}^n \beta_j^2(T)=  \beta_1(S)\beta_1(T)=\frac{\vol(S)\vol(T)}{\vol(V)}.
\end{equation}
This completes the proof of the lemma.
\end{proof}

\begin{lemma}[Expander Mixing Lemma]
\label{lem:expander_mixing}
For any two disjoint subsets $S,T \subseteq V$ with $S\cup T =V$ we have
\begin{equation}
\abs{E(S,T)-\frac{\vol(S)\vol(T)}{\vol(V)}} \leq \lambda \frac{\vol(S)\vol(T)}{\vol(V)},
\end{equation}
where $E(S,T)$ denotes the number edges with one endpoint in $S$ and another in $T$.
\end{lemma}

\begin{proof}
We can write $E(S,T)=\langle \mf{1}(S),\mf{\bar M}\mf{1}(T) \rangle_{\mf D}$. We now write the vectors $\mf{1}(S)$ and $\mf{1}(V)$ in terms of the basis vectors $\mf f_1,\ldots,\mf f_n$ as follows
\begin{align}
    &\mathbf{1}(S) =  \sum_{j=1}^n \beta_j(S) \mf f_j,\\
    &\mathbf{1}(T) = \sum_{j=1}^n \beta_j(T) \mf f_j.
\end{align}
%
Hence, $E(S,T)=\langle \mf{1}(S),\mf{\bar M}\mf{1}(T) \rangle_{\mf D}=\sum_{j=1}^n \beta_j(S)\beta_j(T)\lambda_j$. Furthermore,  by Lemma~\ref{lem:indicator} it follows that
\begin{align}
    E(S,T) = \frac{\vol(S)\vol(T)}{\vol(V)} + \sum_{j=2}^n \beta_j(S) \beta_j(T) \lambda_j.
\end{align}
Hence,
\begin{align}
    \left| E(S,T) -\frac{\vol(S)\vol(T)}{\vol(V)} \right| &=  \left|\sum_{j=2}^n \beta_j(S) \beta_j(T) \lambda_j \right|\\
    & \leq \lambda  \sum_{j=2}^n \abs{\beta_j(S) \beta_j(T) }\label{eq:tmp1}\\
    &=\lambda\sum_{j=2}^n \beta_j^2(S)\\
    &=\lambda \frac{\vol(S)\vol(T)}{\vol(V)},
\end{align}
where the last two lines again follow from Lemma~\ref{lem:indicator}. This completes the proof of the lemma.
\end{proof}

\section{Proof of Lemma~\ref{lem:coupling_bound}}
\label{proof:coupling_bound}

    First, note that in order to generate transitions for the process $\mf X$, we can first generate the transitions for $A(\mf X)$ according to $q_+$ and $q_-$, and then sample the node which will be making the corresponding transition with appropriate probability (e.g, a node $v\in B(\mf x)$ will be sampled with probability $q(\mf x, \mf x+\mf e_{v})/q_+(A(\mf x))$ when the process $A(\mf X(t))$ makes a forward transition). Hence, to prove the lemma, it is sufficient to couple the processes $A(\mf X(t))$ and $\bar A(t)$ such that the desired properties are satisfied.
   We define the coupling as follows. If $\bar A(t) > A(\mf X(t))$, then we let the chains evolve independently. If $\bar A(t) = A(\mf X(t)) = a$ for some $t$, then we construct a coupling to ensure that in the next transition instance $t'$ we have $\bar A(t') \geq A(\mf X(t'))$. In particular, at a common state $a$, in order to generate the forward transition instants for both $A(\mf X)$ and $\bar A$, we first sample two exponential random variables $Z_1$ and $Z_2$ as
\begin{align}
    Z_1 \sim \exp(q_+(a)) \hspace{2em} Z_2 \sim \exp(\bar q_+(\bar a) - q_+(a)),
\end{align}
We then set the next potential forward transition instants of $A(\mf X)$ and $\bar A$ as $t+Z_3$ and $t+Z_4$, where $Z_3=Z_1$ and $Z_4=Z_1\vee Z_2$. 
Similarly, to generate the backward transitions, we first sample $\bar Z_1\sim\exp(\ubar q_-(a))$ and $\bar Z_2\sim\exp(q_-(a)-\ubar q_-(a))$, and then set the potential backward transition instants of $\bar A$ and $A(\mf X)$ as $t+\bar Z_3$ and $t+\bar Z_4$, respectively, where $\bar Z_3=\bar Z_1$ and $\bar Z_4=\bar Z_1 \vee \bar Z_2$. The actual transitions  that occur are decided by events that occur the earliest of the potential events generated above. This coupling ensures that $\bar A(t') \geq A(\mf X(t'))$ at the next transition instant $t'$ and therefore completes the proof. \qed

\section{Proof of Lemma~\ref{lem:f_gen}}
\label{proof:f_gen}

    The fact $\bar q_+(a)-\ubar q_-(a)=\bar q_-(n-a)-\ubar q_+(n-a)$ follows directly from the definitions~\eqref{eq:ubar_plus}-\eqref{eq:ubar_minus}.
    Next, we note that
    \begin{align}
        \bar q_+(a) - \ubar q_-(a)&=
          \frac{1}{2}\,\alpha\bigl(n - a\bigr)
           + (1-\alpha)L_{\max}\,\frac{a^2\bigl(n - a\bigr)}{n^2} - \frac{1}{2}\,\alpha a
           - (1-\alpha)L_{\min}\,\frac{a\bigl(n - a\bigr)^2}{n^2}\\
           & = \frac{1}{2}\alpha (n-2 a) + (1-\alpha)\frac{a\bigl(n - a\bigr)}{n^2} \left [ L_{\max} a - L_{\min}(n-a)\right ]\\
           & = \frac{1}{2}\alpha (n-2 a) - (1-\alpha) \frac{a\bigl(n - a\bigr)}{n^2} \frac{\Sigma_L}{2}(n - 2a) + (1-\alpha) \frac{a\bigl(n - a\bigr)}{n^2} \frac{\Delta_L}{2}n\\
           &\leq (1-\alpha)\frac{\Sigma_L}{2}(n-2a)\brac{\frac{\alpha}{\Sigma_L(1-\alpha)}-\frac{a(n-a)}{n^2}}+ (1-\alpha) \frac{1}{4} \frac{\Delta_L}{2}n\\
           &=n(1-\alpha)\frac{\Sigma_L}{2}F_{\alpha,L}\brac{\frac{a}{n}}.
    \end{align}
    This establishes the first statement of the lemma.
    
    To establish the second part, we take the first derivative of  $F_{\alpha, L}$ given by
    \begin{align}
        F'_{\alpha, L}(y) = -6y^2 + 6y -\left( 1 + \frac{2\alpha}{\Sigma_L (1-\alpha)}\right)=-6\brac{y^2-y+\frac{1}{6}\left( 1 + \frac{2\alpha}{\Sigma_L (1-\alpha)}\right)},
    \end{align}
     which has exactly two real roots under the condition $\alpha < \frac{\Sigma_L}{4+\Sigma_L}$, implied by the condition $\alpha < \frac{\Sigma_L(1-K_L^{1/3})}{4+\Sigma_L(1-K_L^{1/3})}$ stated in the lemma. We set $y_{\min}$ as the minimum of the two roots of the quadratic expression, i.e.
     \begin{equation}
         y_{\min} = \frac{1}{2} - \frac{1}{2\sqrt{3}}\sqrt{1 - \frac{4\alpha}{\Sigma_L(1-\alpha)}}.\label{eq:ymin}
     \end{equation}
     Clearly, $y_{\min} \in (0,1/2)$. It is easy to verify that 
     \begin{align}
         F_{\alpha,L}(y_{\min}) 
         & = \frac{\Delta_L}{4\Sigma_L} - \frac{1}{6\sqrt{3}}\left( 1-\frac{4\alpha}{\Sigma_L(1-\alpha)}\right)^{3/2}.
     \end{align}
     Hence, the condition $F_{\alpha,L}(y_{\min}) < 0$ simplifies to 
     \begin{align}
         \alpha < \frac{\Sigma_L(1-K_L^{1/3})}{4+\Sigma_L(1-K_L^{1/3})},
     \end{align}
which holds under the conditions stated in the lemma. This implies that $F(y_{\min}) < 0$. Furthermore, from the definition of $F_{\alpha,L}$, it is easy to verify that $F_{\alpha,L}(0) >0$ and $F_{\alpha,L}(1/2) > 0$. Since $F_{\alpha,L}$ is a cubic polynomial, the above implies that it has exactly one root $r_1$ in $(0,y_{\min})$ and another root $r_2$ in $(y_{\min},1/2)$. This completes the proof of the first statement of the lemma.

To prove the second statement, we note that for $y > 1/2$ we have $1-2y < 0$. Furthermore, since 
$y(1-y) \leq 1/4$ for all $y\in (0,1)$ and $\alpha > \Sigma_L/(4+\Sigma_L)$, we have $y^2-y+\frac{\alpha}{\Sigma_L(1-\alpha)}\geq \frac{\alpha}{\Sigma_L(1-\alpha)}-\frac{1}{4} > 0$. Hence, we have
\begin{align}
    F_{\alpha,L}(y)& \leq (1-2y)\brac{\frac{\alpha}{\Sigma_L(1-\alpha)}-\frac{1}{4}}+\frac{\Delta_L}{4\Sigma_L},
\end{align}
which simplifies to the RHS of~\eqref{eq:gen_drift_bound}.

\section{Proof for the $2k$-choices model}
\label{proof:2k}

The finite birth and death transition rates for the $2k$-choice rule are
\begin{align}
    &q_+(a) = n \left(\frac{n-a}{n}\right)\left[ \frac{\alpha}{2} + (1-\alpha) \sum_{r=k+1}^{2k} \binom{2k}{r} \left(\frac{a}{n-1}\right)^r \left(\frac{n-a-1}{n-1}\right)^{2k-r}\right],\\
    &q_-(a) = n\left(\frac{a}{n}\right) \left[ \frac{\alpha}{2} + (1-\alpha) \sum_{r=k+1}^{2k} \binom{2k}{r} \left(\frac{n-a}{n-1}\right)^r \left(\frac{a-1}{n-1}\right)^{2k-r}\right].
\end{align}
Notice that $q_-(q) = q_+(n-a)$ and $q_+(a) = q_-(n-a)$. Now define
\begin{align}
    f_{2k,n,\alpha}\left( \frac{a}{n} \right) = \frac{1}{n}(q_+(a)-q_-(a)).
\end{align}
Also define 
\begin{align}
    f_{2k, \alpha}(y) = \lim_{n \rightarrow \infty}f_{2k, n, \alpha}\left( \frac{a}{n} \right),
\end{align}
where
\begin{align}
    f_{2k, \alpha}(y) =  \frac{\alpha}{2}(1-2y) + (1-\alpha) \sum_{r=k+1}^{2k}\binom{2k}{r}\left[y^r (1-y)^{2k-r+1} - (1-y)^r y^{2k-r+1}\right].
\end{align}
Clearly, $y=1/2$ is a root of $f_{2k, \alpha}(y)$. We also have $f_{2k, \alpha}(0) = \alpha/2$ and $f_{2k, \alpha}(1)=-\alpha/2$. We note here that $f_{2k, \alpha}$ is continuous as it is a polynomial. We can write $f_{2k, \alpha}$ alternatively in terms of incomplete beta functions, i.e. 
\begin{align}
    f_{2k, \alpha}(y) = \frac{1}{2}\alpha (1-2y) + (1-\alpha)[(1-y) I_y(k+1,k) - yI_{1-y}(k, k+1)],
\end{align}
where 
\begin{align}
    I_y(a,b) & \overset{(c)}{=} \sum_{r=a}^{a+b-1} \binom{a+b-1}{r}y^r(1-y)^{a+b-1-r}\\
    & \overset{(d)}{=} \frac{1}{B(a,b)}\int_0^y t^{a-1}(1-t)^{b-1}dt,
\end{align}
is the incomplete beta function with shape parameters $a$ and $b$, and $B(a,b)$ is the complete beta function with shape parameters $a$ and $b$. We will use both forms of equality, (c) and (d), in this section. Using the fact that $I_{1-y}(k+1, k) = 1 - I_{y}(k, k+1)$, we can write
\begin{align}
    f_{2k, \alpha}(y) = \frac{1}{2}\alpha (1-2y) + (1-\alpha)[-y + (1-y) I_y(k+1,k) + yI_{1-y}(k, k+1)].
\end{align}
For this analysis, we require the first two derivatives of this function. The first derivative of the function can be written as
\begin{align}
\label{2k function derivative}
    f_{2k, \alpha}'(y) = -\alpha + (1-\alpha)\bigg [-1  + I_y(k,k+1) &- I_x(k+1,k)\\
    & + (1-y)\frac{dI_y(k+1,k)}{dy}\\
    & + y\frac{dI_y(k,k+1)}{dy} \bigg ].
\end{align}
We now write
\begin{align}
    \frac{dI_y(k+1,k)}{dy} &= \frac{1}{B(k+1,k)}\frac{d}{dy}\int_0^{y}t^{k}(1-t)^{k-1}dt\\
    & = \frac{y^k(1-y)^{k-1}}{B(k+1, k)}.
\end{align}
We can now use the following
\begin{align}
    B(k+1, k) &= \frac{\Gamma(k+1)\Gamma(k)}{\Gamma(2k+1)}\\
    & = \frac{k! (k-1)!}{(2k)!}\\
    & = \frac{1}{k \binom{2k}{k}}.
\end{align}
Hence, we can write
\begin{align}
    \frac{dI_y(k+1,k)}{dy} &= k\binom{2k}{k}y^k(1-y)^{k-1},
\end{align}
and similarly for 
\begin{align}
    \frac{dI_y(k,k+1)}{dy} = k \binom{2k}{k}y^{k-1}(1-y)^{k}.
\end{align}
Plugging these into equation \eqref{2k function derivative}, we obtain
\begin{align}
    f_{2k, \alpha}'(x) = -\alpha + (1-\alpha)\left \{-1  + I_y(k,k+1) - I_y(k+1,k) + 2k\binom{2k}{k}y^k(1-y)^k \right \}.
\end{align}
Now we focus on $I_y(k,k+1) - I_y(k+1,k)$, which can be written as
\begin{align}
    I_y(k,k+1) - I_y(k+1,k) &= \frac{1}{B(k, k+1)}\left[\int_{0}^{y} t^{k-1}(1-t)^k dt - \int_{0}^{x} t^{k}(1-t)^{k-1} dt \right]\\
    & = \frac{1}{B(k, k+1)}\left[\int_{0}^{y} (t^{k-1}(1-t)^k - t^{k}(1-t)^{k-1} )dt \right]\\
    & = \frac{1}{B(k, k+1)}\left[\int_{0}^{y} t^{k-1}(1-t)^{k-1}(1-2t)dt \right]\\
    & = \frac{1}{B(k, k+1)}\frac{1}{k}y^{k}(1-y)^{k}\\
    & = \binom{2k}{k}y^{k}(1-y)^{k}.
\end{align}
Hence, plugging this into equation \eqref{2k function derivative} and simplifying, we finally end up with
\begin{align}
    f_{2k, \alpha}'(y) -1 + (1-\alpha)(2k+1)\binom{2k}{k}y^{k}(1-y)^{k}.
\end{align}
The second derivative follows easily from here, i.e.
\begin{align}
    &f_{2k, \alpha}''(y) = (1-\alpha)k(2k+1)\binom{2k}{k}(1-2y)y^{k-1}(1-y)^{k-1}.
\end{align}
Now we prove some statements about $f'_{2k, \alpha}(y)$. It is easy to see that the derivative is even about $y=1/2$, i.e. $f'_{2k, \alpha}(y) = f'_{2k, \alpha}(1-y)$. The derivative evaluated at both $y=0$ and $y=1$ is $f'_{2k, \alpha}(0) = f'_{2k, \alpha}(1) = -1$. The derivative also has a local optima when $f''_{2k, \alpha}(y)=0$, i.e. when $y=0, 1/2, 1$. It is obvious that for every $y \in (0, 1/2)$ that $f''_{2k, \alpha}(y) > 0$, and hence, $f'_{2k, \alpha}(y)$ is strictly increasing in this region, and therefore by the even property, $f'_{2k, \alpha}(y)$ is decreasing on $(1/2,1)$. Therefore, we conclude that $f'_{2k, \alpha}(1/2)$ is a local optimum, and the curve $f'_{2k, \alpha}(y)$ is unimodal on $[0,1]$. The value of this optimum is 
\begin{align}
    f'_{2k, \alpha}(1/2) = (1-\alpha)(2k+1)\binom{2k}{k}\left( \frac{1}{2}\right)^{2k} - 1.
\end{align}
Now let 
\begin{align}
    \alpha_{2k} = 1- \frac{2^{2k}}{(2k+1) \binom{2k}{k}}.
\end{align}

\begin{lemma}
    The function $f_{2k, \alpha}$ satisfies the following properties
    \begin{enumerate}
        \item  $f_{2k,\alpha}\left(y\right) = - f_{2k, \alpha}(1-y)$.
        \item Fix $\alpha < \alpha_{2k}$. Then there exists a region $(r_-, 1/2)$ such that $f_{2k,\alpha}(y) < 0$ when $y\in (r_{-},1/2)$. 
        \item For all $\alpha \in (0,1)$, we have $f_{2k,\alpha}(y)-f_{n,\alpha}(y')\leq -c(2k,\alpha) (y-y')$ for any $0\leq y'\leq y\leq 1$, where $c(2k,\alpha) > 0$ for all $k \geq 1$ and $\alpha > \alpha_{2k}$.  
\end{enumerate}
\end{lemma}

\begin{proof}
    \textbf{Statement 1:} The result in statement 1 actually holds for any finite $n$. It follows directly from the symmetry of the transition rates. We can write
    \begin{align}
        f_{2k, n, \alpha}\left(1-\frac{a}{n}\right) &= f_{2k, n, \alpha}\left(\frac{n-a}{n}\right)\\
        & = \frac{1}{n}(q_+(n-a) - q_-(n-a))\\
        & \overset{(a)}{=} \frac{1}{n}(q_-(a) - q_+(a))\\
        & = - \frac{1}{n}(q_+(a) - q_-(a))\\
        & = - f_{2k, n, \alpha}\left(\frac{a}{n}\right),
    \end{align}
    where $(a)$ has used $q_-(q) = q_+(n-a)$ and $q_+(a) = q_-(n-a)$. This concludes the proof of the first statement of the lemma. 

    \textbf{Statement 2:} The local optimum of $f'_{2k, \alpha}$ is greater than zero when $\alpha < \alpha_{2k}$, and we know that $f'_{2k, \alpha}(0)=f'_{2k, \alpha}(1)=-1$, as well as the fact that $f'_{2k, \alpha}(y)$ is unimodal on [0,1]. Hence there are exactly 2 roots of $f'_{2k, \alpha}(y)$ in [0,1]. Call these roots $y_-$ and $y_+$, with $y_- = 1-y_+$ due to the evenness of $f'_{2k, \alpha}$ around $y=1/2$. So $f_{2k,\alpha}'(y) < 0$ for all $y \in (0, y_-) \cup(y_+,1)$, and $f_{2k, \alpha}'(y) > 0$ for all $y \in (y_-, y_+)$. Hence, $f_{2k, \alpha}(y)$ has a local minimum and maximum at $y_-$ and $y_+$ respectively. Since $f_{2k, \alpha}(y)$ is increasing on $(y_-,1/2)$, we have $f_{2k, \alpha}(y_-) < f_{2k, \alpha}(1/2) = 0$. Since $f_{2k, \alpha}(0) > 0$ and $f_{2k, \alpha}(y_-) < 0$, the Intermediate Value Theorem gives a root of $f_{2k, \alpha}$ at $r_- \in (0, y_-)$. Using the first statement of this lemma, we can immediately conclude the existence of another root at $r_+ = 1-r_- \in (y_+,1)$. Together with the root at $y=1/2$, we deduce $f_{2k, \alpha}(y)$ has exactly 3 roots on [0,1] for $\alpha < \alpha_{2k}$. Therefore, $f_{2k, \alpha}(y)<0$ for $y \in (r_-, 1/2)$. This concludes the proof of the second statement. 

    \textbf{Statement 3:} The local optimum of $f'_{2k, \alpha}$ is  negative when $\alpha > \alpha_{2k}$, and we know that $f'_{2k, \alpha}(0)=f'_{2k, \alpha}(1)=-1$, as well as the fact that $f'_{2k, \alpha}(y)$ is unimodal on [0,1]. Therefore in this regime, $f'_{2k, \alpha}(y) < 0$ for every $y \in[0,1]$. Hence, there exists a constant such that $f'_{2k, \alpha}(y) < -c(2k, \alpha)$ for all $y \in (0,1)$. Then we can say that for $y, y' \in (0,1)$, with $y>y'$, we have 
    \begin{align}
        f_{2k, \alpha}(y) - f_{2k, \alpha}(y') &= \int_{y'}^{y} f'_{2k, \alpha}(s)ds\\
        & \leq -c(2k, \alpha)\int_{y'}^{y}ds\\
        & = -c(2k, \alpha)(y-y').
    \end{align}
    This concludes the proof of the third statement of the lemma.
\end{proof}

\end{document}